 \documentclass[10pt]{article}
 \usepackage[top=1in,bottom=1in,left=1.5in,right=1.5in]{geometry}
  \usepackage{amsmath,amssymb,amsthm}
  \usepackage{latexsym}% defines $\Box$ for LaTeX2e
  \usepackage[dvips]{pstcol} % PSTricks
  \usepackage[dvips]{graphicx}
  \allowdisplaybreaks
  \newcommand{\C}{\mathbb{C}}
  \newcommand{\F}{\mathbb{F}}
  \newcommand{\N}{\mathbb{N}}
  \renewcommand{\P}{\mathbb{P}}
  \newcommand{\R}{\mathbb{R}}
  \newcommand{\Z}{\mathbb{Z}}

  \newcommand{\h}{\mathbf{H}}

  \newcommand{\x}{\mathbf{x}}

  \newcommand{\0}{\mathbf{0}}
  \newcommand{\1}{\mathbf{1}}
  
  \newcommand{\Gr}     {{\rm Gr\;}}

  \newcommand{\rM}{\mathrm{M}}

  \newcommand{\rA}{\mathrm{A}}

  \newcommand{\rS}{\mathrm{S}}
  
  \newcommand{\rU}{\mathrm{U}}
  \newcommand{\rV}{\mathrm{V}}
  \newcommand{\rW}{\mathrm{W}}

  \newcommand{\hs}{\hspace*{\parindent}}
  \newcommand{\cl}{\mathop{\mathrm{Cl}}\nolimits}
  \def\gl{{\rm GL}}

  \newcommand{\trans}{^\top}

  \def\rank{\mathop{{\rm rank\;}}\nolimits}
  \newcommand{\codim}{\mathrm{codim\;}}

  \newtheorem{theo}{\bfseries \hs Theorem}[section]
  
  \newtheorem{prop}[theo]{\bfseries \hs Proposition}
  \newtheorem{lemma}[theo]{\bfseries \hs Lemma}
  \newtheorem{corol}[theo]{\bfseries \hs Corollary}
  
  \newtheorem{problem}[theo]{\bfseries \hs Problem}

 \theoremstyle{remark}
  \newtheorem{rem}[theo]{\rm\bfseries \hs Remark}
  \newtheorem{example}[theo]{\rm\bfseries \hs Example}

 \numberwithin{equation}{section} % Automatically number equations within sections

 \def\({\left(}

 \def\){\right)}

 \def\cl#1{\left\lceil#1\right\rceil}

 \def\th{\theta}

\begin{document}

 \title{$2$-adic
 valuations of certain
 ratios of products \\ of
 factorials and applications}

 \author
 {Shmuel Friedland$^1$ and Christian Krattenthaler$^{2}$%
 \thanks{Research partially supported by
 EC's IHRP Programme, grant HPRN-CT-2001-00272,
``Algebraic Combinatorics in Europe"\newline
 \indent\kern3pt$^\dagger$Current address: Fakult\"at f\"ur Mathematik, Universit\"at Wien,
 Nordbergstra{\ss}e~15, A-1090 Vienna, Austria} $^\dagger$
 \\\\
 $^1$\it Department of Mathematics, Statistics and Computer Science,\\
 \it University of Illinois at Chicago \\
 \it Chicago, Illinois 60607-7045, USA\\
  Email: {\tt friedlan@uic.edu}\\
 \\
 $^2$\it Institut Camille Jordan,
 \it Universit\'e Claude Bernard Lyon-I\\
 \it 21, avenue Claude Bernard,\\\it F-69622 Villeurbanne Cedex,
 France\\
 \centerline{WWW: \tt
 http://igd.univ-lyon1.fr/\~{}kratt}
 }

 \date{August 24, 2005}
 \maketitle

 \begin{abstract}
 We prove
 the conjecture of Falikman--Friedland--Loewy
 on the parity of the degrees of projective varieties of $n\times
 n$ complex symmetric matrices of
 rank at most $k$.
 We also characterize the parity of the degrees of projective varieties of $n\times
 n$ complex skew symmetric matrices of
 rank at most $2p$.
 We give recursive relations which determine the parity of the
 degrees of projective varieties of $m\times n$ complex matrices
 of
 rank at most $k$.
 In the case the degrees of these varieties are odd, we characterize
 the minimal dimensions of subspaces of $n\times
 n$ skew symmetric real matrices and of $m\times n$
 real matrices containing a nonzero matrix of rank
 at most $k$.
 The parity questions studied here are also of combinatorial interest
 since they concern the parity of the number of plane partitions
 contained in a given box, on the one hand, and the parity of
 the number of symplectic tableaux of rectangular shape, on the
 other hand.
 \end{abstract}

\bigskip
 \noindent {\bf 2000
 Mathematics Subject Classification:}
{\it Primary}: 11A07, 15A03; \\{\it Secondary}: 14M12, 14P25,  15A30.
 \par\noindent {\bf Keywords:}  2-adic
 valuations of ratio of products
 of factorials, parity of degrees of determinantal varieties,
 subspaces of real skew symmetric matrices, subspaces of real rectangular matrices,
 parity of number of plane partitions, parity of number of symplectic
 tableaux.

 \section{Introduction}
  Consider the polynomial equation $z^{d}+a_1z^{d-1}+\dots+a_d=0$
  over the field of complex numbers $\C$.  The fundamental theorem
  of algebra
  says that this polynomial system has always a
  nontrivial complex solution $\zeta$.  Assume that $a_1,\ldots,a_d$ are
  real valued.  Clearly, this does not imply that the polynomial
  equation is solvable over the field of real numbers $\R$.
  Since the complex solutions come in complex pairs,
  it is well-known that
  a sufficient condition for an existence
  of a real solution $\zeta$ is that the degree $d$ is odd.
  A similar situation holds in
  a more general setting.

  Let $\P\R^n$ and $\P^n$ be the real and
  the complex projective space of dimension $n$, respectively. Let $V
  \subset \P\R^n$ be an algebraic variety, such that
  its
  complexification $V_{\C} \subset \P^n $ is irreducible and has
  codimension $m \ge 1$.  Hence $B(V_{\C})=V_{\C}$, where $B:\P^n \to \P^n$ is
  the involution $z \mapsto \bar z$.  Recall that any linear space
  $M \subset \P^n$ of dimension $m$ intersects $V_{\C}$. Furthermore, for
  a generic $M$, the set $V_{\C} \cap M$ consists of
  a finite number, $d$ say, of points.
  The positive integer $d$ is called the degree of
  $V_{\C}$ and is denoted by $\deg V_{\C}$. As in the previous
  simple case, it is well-known that if $\deg V_{\C}$ is odd then any linear
  space $L \subset \P\R^n$ of dimension $m$ intersects $V$. Indeed,
  for a generic $L \subset \P\R^n$, the set $V_{\C}\cap L_{\C}$
  consists of $\deg V_{\C}$ points.
  As this set is invariant under the involution $B$,
  we deduce that there exists $\zeta \in V_{\C} \cap L_{\C}$ such
  that $B(\zeta)=\zeta$,
  that is, $\zeta \in \P\R^n$.
  A continuity
  argument yields that $V\cap L \ne \emptyset$ for any $L$ in the
  Gra\ss mannian $\Gr(m+1,n+1,\R)$.  We refer to this result as the
  \emph{odd degree theorem}.

 Assume now that $\deg V_{\C}$ is even.  Then it is not difficult
 to find nontrivial examples where $V \cap L' =\emptyset$ for some
 $L'\in \Gr(m+1,n+1,\R)$. Moreover examples can be found among
 determinantal varieties, such that for any integer $k \in [0,p)$
 there exists $L' \in \Gr(m+k+1,n+1,\R)$ satisfying $V\cap L'
 =\emptyset$, while $V\cap L \ne\emptyset$ for any $L \in
 \Gr(m+p+1,n+1,\R)$ (see \cite[\S2]{FL}).

 The following generalization of the odd degree theorem is
 proved
 in \cite{FL}: Assume that $\deg V_{\C}$ is even and let $r$ be a
 positive integer. Suppose that the codimension of the variety of
 the singular points of $V_{\C}$ in $V_{\C}$ is at least $2r+1$.
 Suppose furthermore that for a generic $M \in \Gr(m+2r+1,n+1,\C)$
 the Euler characteristic of $V_{\C}\cap M$ is odd.  Then $V\cap L
 \ne\emptyset$ for any $L \in \Gr(m+2r+1,n+1,\R)$.

 The degree of $V_{\C}$ and the Euler characteristic of $V_{\C}\cap
 M$ can be often computed using the projectivized complex bundles
 and the corresponding Chern classes of their tangent bundles, see
 for example \cite{HT2} and \cite{FL}. It turns out that the
 degrees of determinantal varieties (see
 the special cases
 discussed below) are quotients of products of certain binomial
 coefficients (see \cite{HT1} and \cite[Ch. II]{ACGH}). Binomial
 coefficients appear also in
 Stiefel--Whitney classes \cite{MS}.

 Let $\rM_{m,n}(\F)$, $\rM_n(\F)$, $\rS_n(\F)$, $\rA_n(\F)$ be the spaces of
 $m\times n$ matrices,
 $n\times n$ matrices, $n \times n$ symmetric matrices,
 and $n\times n$ skew (antisymmetric) matrices
 with entries in $\F=\R,\C$, respectively. Let
 $\rU_{k,m,n}(\F)$, $\rV_{k,n}(\F)$, $\rW_{k,n}(\F)$ be the
 varieties of all matrices in
 $\rM_{m,n}(\F)$, $\rS_n(\F)$, $\rA_n(\F)$ of rank $k \le \min(m,n)$ or less, respectively.
 It is known that any $A\in \rA_n(\F)$ has an even rank,
 see e.g.\
 \cite[\S11.4]{Gan} or \S5.  Hence it is enough to consider
 $\rW_{2p,n}(\F)$, where $1\le p\le \lfloor \frac{n}{2}\rfloor$.
 Then the
 projectivizations $\P\rU_{k,m,n}$, $\P \rV_{k,n}(\F)$, $\P \rW_{2p,n}(\F)$
 are irreducible
 varieties of
 codimension $(m-k)(n-k)$, $\binom {n-k+1} 2$, $\binom {n-2p} 2$
in the projective
 spaces
 $\P\rM_{m,n}(\F)$, $\P \rS_n(\F)$, $\P \rA_n(\F)$, respectively.  Note that
 $\P \rU_{k-1,m,n}(\F)$, $\P \rV_{k-1,n}(\F)$,
 $\P \rW_{2(p-1),n}(\F)$ are the
 varieties of the
 singular points of
 $\P \rU_{k,m,n}(\F)$, $\P \rV_{k,n}(\F)$, $\P \rW_{2p,n}(\F)$, respectively.
 For $\P\rU_{k,m,n}(\F)$, $\P \rV_{k,n}(\F)$ see for example \cite[II]{ACGH},
 and for $\P \rW_{2p,n}(\F)$ see \S5.
 Let $d(m,n,k,\F)$, $d_s(n,k,\F)$, $d_a(n,2p,\F)$ be the smallest integer
 $\ell$ such that every $\ell$-dimensional subspace of $\rM_{m,n}(\F)$, $\rS_n(\F)$, $\rA_n(\F)$
 contains a nonzero matrix
 whose rank is at most $k, 2p$, respectively.  Then
 \begin{equation}\label{dnkc}
 \kern-2pt
 d(m,n,k,\C)=(m-k)(n-k)+1, \;
 d_s(n,k,\C)=\binom {n-k+1} 2 +1, \; d_a(n,2p,\C)=\binom {n-2p} 2 +1,
 \kern-2pt
 \end{equation}
 and the problem is to determine $d(m,n,k,\R)$, $d_s(n,k,\R)$, $d_a(n,2p,\R)$.
 The degrees of
 $\P\rU_{k,m,n}(\C)$, $\P\rV_{k,n}(\C)$, $\P \rW_{2p,n}(\C)$
 were computed by Harris and Tu in \cite{HT1},
 \begin{eqnarray}
 &&\gamma_{k,m,n}:=\deg
 \P\rU_{k,m,n}(\C)=
 \prod _{j=0} ^{n-k-1}\frac {\binom {m+j}{m-k}} {\binom {m-k+j}{m-k}}=
 \prod_{j=0}^{n-k-1}\frac{(m+j)!\,j!}{(k+j)!\,(m-k+j)!},
 \label{degkmn}\\
 &&\delta_{k,n}:={\rm deg\;} \P \rV_{k,n}(\C)=\prod_{j=0}^{n-k-1} \frac{\binom {n+j}
 {n-k-j}}{\binom {2j+1} j}, \quad \varepsilon_{2p,n}:= \deg \P
 \rW_{2p,n}=\frac{\delta_{2p+1,n}}{2^{n-2p-1}}.
 \label{degvk}
 \end{eqnarray}
 For the curiosity of the reader we remark that these quantities
 have also combinatorial interpretations. The quantity $\gamma_{k,m,n}$
 counts plane partitions which are contained in an
 $(n-k)\times(m-k)\times k$ box (see \cite{Bressoud} and \S6). On the other
 hand, the quantity $\varepsilon_{2p,n}$ counts symplectic tableaux
 (see \cite{King}) of a rectangular shape of size $n\times p$, and
 thus several other sets of combinatorial objects
 (see \cite{Sundaram} and \cite{GhorKrat}
 for more information on these topics).

 It was shown in \cite{FFL} that $\delta_{n-q,n}$ is odd if
 \begin{equation}\label{spc1}
 n\equiv \pm q \;({\rm mod\;}2^{\lceil \log_2 2q\rceil}) .
 \end{equation}
 For values of $q$ and $n$ which satisfy this condition,
 \begin{equation}\label{spc}
 d_s(n,n-q,\C)=d_s(n,n-q,\R)= \binom {q+1} 2 +1.
 \end{equation}
 It was
 furthermore shown in \cite{FFL} that this equality implies the following
 interesting result.  Assume that $n \ge q$ and
 that $n$ satisfies
 (\ref{spc1}), then any
 $\binom {q+1} 2$-dimensional subspace of $\rS_n(\R)$
 contains a nonzero matrix with an eigenvalue of multiplicity at
 least $q$. This statement for $q=2$ yields
 Lax's
 result \cite{Lax} that any $3$-dimensional subspace of $\rS_n(\R)$
 contains a nonzero matrix with a multiple eigenvalue for $n\equiv
 2 \;({\rm mod\;} 4)$.  (This result and its generalization in
 \cite{FRS} is of importance in the study of singularities of
 hyperbolic systems.)

 On the other hand,
 it was conjectured in \cite{FFL} that
 also the converse holds, that is, that if $\delta_{n-q,n}$ is odd then
 (\ref{spc1}) holds.  In particular, for $n$ and $q$ which do not
 satisfy (\ref{spc1}), we do not have a \emph{simple} way to
 compute $d_s(n,n-q,\R)$.
 One of the main purposes of this
 paper is to prove this conjecture,
 see \S\S3,4.  Our results yield
 in addition that
 $\varepsilon_{2p,n}$ is odd if and only
 if (\ref{spc1}) holds with $q=n-2p$,
 see \S5. Hence for these values of $p,n$
 we have
 \begin{equation}\label{asrc}
 d_a(n,2p,\C)=d_a(n,2p,\R)= \binom {n-2p} 2 +1.
 \end{equation}
 In particular, for $n\equiv
 2 \;({\rm mod\;} 4)$ any two-dimensional subspace of real
 $n\times n$ skew symmetric matrices contains a nonzero singular
 matrix.
 For $n$ and $q=n-2p$ which do not satisfy the condition
 (\ref{spc1}), we do not have a \emph{simple} way to
 compute $d_a(n,2p,\R)$.

 We also consider the problem
 of characterizing the values of $k,m,n$ for which
 $\gamma_{k,m,n}$ is odd.
 This problem seems to have a rather intricate solution.
 We give some partial results on this problem
 in \S6. In particular, we provide an algorithm for computing the
 parity of $\gamma_{k,m,n}$ from the binary expansions of $k,m,n$
 directly, without having to actually compute $\gamma_{k,m,n}$
 (see Remark~\ref{rem2} and Proposition~\ref{oddconb}). As above,
 if $\gamma_{k,m,n}$ is odd then
 \begin{equation}\label{dmnkcr}
 d(m,n,k,\C)=d(m,n,k,\R)=(m-k)(n-k)+1.
 \end{equation}
 See Corollary~\ref{cordkmn} for the corresponding geometric results
 which we are able to prove.

 Another purpose of this paper
 is to estimate the
 $2$-adic valuation of $\delta_{k,n}$, i.e., the
 largest power of $2$ that divides $\delta_{k,n}$.
 There are two reasons for these estimations.
 First, we show that our methods can give good estimates
 for the complex behavior of the $2$-adic valuation of
 $\delta_{k,n}$.  Second, we recall the classical results
 of Radon \cite{Rad} and Hurwitz \cite{Hur} on the maximal dimension of the spaces of
 $n\times n$ scaled orthogonal matrices,
 and the Adams result \cite{Ada} on maximal number of vectors fields on the tangent bundle of
 $n-1$-dimensional sphere,
 which are functions of the $2$-adic valuation of $n$.
 We believe that the $2$-adic
 valuation of $\deg \P V_{k,n}(\C)$
 is related to a lower bound for the problem raised in
 Friedland--Libgober \cite{FL}.

 Consider the variety of $n\times n$ singular
 matrices in $\rM_n(\F)$.  Clearly the degree
 of this variety is
 $n$ and
 its codimension is $1$.  Hence any two-dimensional complex
 subspace $L\subset \rM_n(\C)$ contains a nonzero singular matrix.
 For $n$ odd, any two-dimensional real subspace of $\rM_n(\R)$
 contains a nonzero singular matrix.  For $n$ even, the situation
 is much more complicated.
 For $n\in\N$, let $c+4d$ be the $2$-adic valuation of $n$. Thus
 $n=(2a+1)2^{c+4d}$, where $a$ and $d$ are nonnegative integers and
 $c\in\{0,1,2,3\}$. Then the Radon--Hurwitz number $\rho(n)$ is
 defined by $\rho(n)=2^c+ 8d$.
 The classical results of Radon \cite{Rad} and Hurwitz \cite{Hur}
 state that $\rho(n)$ is the maximal dimension of an
 $n$-dimensional subspace
 $U$ of $\rM_n(\R)$
 such that each
 nonzero $A \in U$ is an orthogonal
 matrix times $r \in \R^*$. In his famous work \cite{Ada}, Adams
 gave a nonlinear version of the Radon--Hurwitz
 result by showing
 that $\rho(n)-1$ is the maximal number of linearly independent
 vector fields on the $(n-1)$-dimensional sphere in $\R^n$.
 In particular,
 Adams's result implies that any $(\rho(n)+1)$-dimensional subspace of
 $\rM_n(\R)$ contains a nonzero singular matrix.

 The parity of binomial coefficients plays
 a role in generalized Radon--Hurwitz
 numbers \cite[Prop. I (f)]{BF}.
 Similarly, we believe that the answer to the following problem
 raised in \cite{FL} significantly depends on
 the $2$-adic valuation of $\delta_{k,n}$:

 \begin{problem}\label{mprob}  Assume that $\delta_{k,n}$ is even.  Find
 an integer $r \ge 1$, preferably the smallest possible, such that
 $ 2r < \binom {n-k+2} 2 - \binom {n-k+1} 2$,
 and
 such that the Euler characteristic of $\P \rV_{n,k}(\C) \cap M$ is odd for
 a generic $M \in \Gr\!\left(\binom {n-k+1}
 2+2r+1,\binom {n+1} 2,\C\right)$.
 \end{problem}
 For the above minimal value of $r$,
 we have $\P \rV_{n,k}(\R) \cap
 L\ne\emptyset$ for any
$$L \in \Gr\!\left(\binom {n-k+1}
 2+2r+1,\binom {n+1} 2,\R\right).$$

 We now briefly survey the contents of this paper.  In \S2 we give
 some
 auxiliary results on the 2-adic
 valuation of $\delta_{n-q,n}$.
 In \S3 we prove the conjecture
 from \cite{FFL} characterizing the values of $q$ and $n$ for which
  $\delta_{n-q,n}$ is odd.  In \S4 we discuss
 properties of the 2-adic
 valuation of $\delta_{n-q,n}$ when the condition
 (\ref{spc1}) is not satisfied. In particular, we characterize the
 values of $q$ and $n$ for which the 2-adic valuation of
 $\delta_{n-q,n}$ is $1$, and we show that, for fixed $q$,
 the $2$-adic valuations of $\delta_{n-q,n}$ have a wave-like
 behavior as $n$ increases. In \S5 we show that
 $\varepsilon_{2p,n}$ is odd if and only if
 (\ref{spc1}) holds with $q=n-2p$.
 Hence for these values of $n$ and $p$
 the equality (\ref{asrc}) holds.  Finally, in \S6 we study
 the parity of the number of plane partitions contained in an $a\times
 b\times c$ box, and thus the parity of $\gamma_{k,m,n}$. Some partial
 results are given, as well as the above-mentioned algorithm for
 computing this parity.

 \section{Preliminary results
 on the $2$-adic valuation of $\delta_{n-q,n}$}
 For a nonzero integer $i$ we write $\nu_2(i)$ for the $2$-adic
 valuation of $i$.  That is $i=(2j+1)2^{\nu_2(i)}$ for some
 $j\in\Z$.
 For positive integers $q$ and $n$ define
 \begin{equation}\label{deftqn}
 \th_{q,n}:=\prod _{j=0} ^{q-1}\frac {\binom {n+j} {q-j}} {\binom
 {2j+1}{j}}. %\quad \big(=\delta_{n-q,n}\big).
 \end{equation}
 The reader should note that $\th_{q,n}=\delta_{n-q,n}$
 (compare (\ref{degvk})).
 It will be convenient later to extend the definition of $\th_{q,n}$
 to all {\it nonnegative} integers $q$ and $n$, that is, to allow
 $n=0$, respectively $q=0$, in (\ref{deftqn}) also. In particular, for
 $q=0$ we set $\th_{0,n}=1$ by interpreting the empty product as $1$.

 Clearly, we have $\theta_{q,n}=0$ for $q>n$.
 We want to study the behavior of the $2$-adic
 valuation of
 $\th_{q,n}$ for $n\ge q$. The following proposition simplifies this study
 as it exhibits a simple relationship
 between the $2$-adic
 valuation of $\th_{q,n}$ when $n-q$ is even
 and those when $n-q$ is odd. In particular, this result allows one
 to concentrate on the analysis of just one case, which will be the
 case where $n-q$ is even.

 \begin{prop} \label{prop:1}
 Let $n$ and $q$ be
 nonnegative integers, $n\ge q+1$. Then
 \begin{equation}\label{recrnt}
 \nu_2(\th_{q,n})=\nu_2(\th_{q+1,n})+q -\sum_{j=0}^q\nu_2(n-q+2j).
 \end{equation}
 In particular, if $n-q$ is odd then
 $\nu_2(\th_{q,n})=\nu_2(\th_{q+1,n})+q\ge q$.
 Hence, if $n$ and $q$ are both positive, and if $n-q$ is positive and
 odd, then $\th_{q,n}$ is even.
 \end{prop}
 \begin{proof}
 The ratio of $\th_{q+1,n}$ and $\th_{q,n}$ is
 \begin{align*}
 \frac {\th_{q+1,n}} {\th_{q,n}}&=\frac{1}{\binom {2q+1} q}
 \prod_{j=0}^q \frac{(n+j)!}{(q+1-j)!\,(n-q-1+2j)!}
 \prod_{j=0}^{q-1}\frac{(q-j)!\,(n-q+2j)!}{(n+j)!}\\
 &=\frac{(q)!\,(q+1)!}{(2q+1)!}\frac{(n+q)!}{1!\,(n+q-1)!}\prod_{j=0}^{q-1}
 \frac{(q-j)!\,(n-q+2j)!}{(q+1-j)!\,(n-q-1+2j)!}\\
 &=\frac{(q+1)!}{2^q(2q+1)!!} (n+q)
 \prod_{j=0}^{q-1}
 \frac{(n-q+2j)}{(q+1-j)}=\frac { \prod _{j=0} ^{q}(n-q+2j)}
 {2^q\,(2q+1)!!}.
 \end{align*}
 (Here $(2q+1)!!:=(2q+1)\cdot(2q-1)\cdots 3\cdot 1$.)
 As $(2q+1)!!$ is odd, we deduce (\ref{recrnt}).
 Assume that $n-q$ is odd.  Then
 $n-q+2j$ is odd for $j=0,\ldots,q$, and the last part of the proposition follows.
 \end{proof}

 We now concentrate on the case where $n-q$ is even.

 \begin{prop}\label{eqfthm1}
 Let $n$ and $q$ be
 nonnegative integers, $n\ge q$, such that the difference
 of $n$ and $q$ is even, say $n-q= 2p$.  Then
 \begin{equation}\label{forthqn}
 \nu_2(\th_{q,n})=(n-1-p)p-\nu_2\left(\prod_{k=1}^p\frac{(n-k)!}{(k-1)!}\right).
 \end{equation}
 Here
 again, in case that $p=0$, the empty product has to be interpreted as $1$.
 In particular, $(n-1-p)p\ge \nu_2\left(\prod_{k=1}^p\frac{(n-k)!}{(k-1)!}\right)$.
 \end{prop}
 \begin{proof}

 We prove (\ref{forthqn}) by a reverse induction on $q$. By the
 definition (\ref{deftqn}) of $\th_{q,n}$ we have $\th_{n,n}=1$.
 Hence $\nu_2(\th_{n,n})=0$, which confirms (\ref{forthqn}) in this case.

 Proposition~\ref{prop:1} implies that for any positive integer $k$
 we have
 $$\nu_2(\th_{n-2k+1,n})=\nu_2(\th_{n-2k+2,n})+n-2k+1.$$
 We now use
 (\ref{recrnt}) for $q=n-2k$ to obtain the recursive formula
 \begin{align*}
 \nu_2(\th_{n-2k,n})&=\nu_2(\th_{n-2k+2,n})+ n-2k+1 +n-2k -
 \sum_{j=0}^{n-2k}\nu_2(2k+2j)\\
 &=\nu_2(\th_{n-2k+2,n})+n-2k -\sum_{j=0}^{n-2k}\nu_2(k+j)\\
 &=
 \nu_2(\th_{n-2k+2,n})+n-2k -\nu_2\left(\frac{(n-k)!}{(k-1)!}\right).
 \end{align*}
 Use this recursive relation for $k=p,p-1,\ldots,1$ to obtain
 (\ref{forthqn}).  Since $\nu_2(\theta_{q,n})\ge 0$ we obtain
 that the right-hand side of (\ref{forthqn}) is nonnegative.
\end{proof}

 Our next goal is to give an explicit expression of the $2$-adic
 valuation of $\th_{q,n}$ in terms of binary digit sums. More
 precisely, for a nonnegative integer $m$ let $s(m)$ denote the sum of the digits of $m$
 when written in binary notation.  Then $s(0)=0$, $s(2m)=s(m)$,
 $s(2m+1)=1+s(m)$, and
 \begin{equation}\label{AI}
 s(2^e-1-k)=e-s(k) \quad \text{for any integers } e\ge 0,\; k\in [0,2^e-1].
 \end{equation}
 The basic result which ties together the $2$-adic valuation of
 factorials and binary digit sums is the following one due to Legendre
 (cf.\ \cite[Sec.~4.4]{GrKnPa} and \cite{Legendre}).
 We bring its proof for completeness.

 \begin{prop}\label{2adn!}  Let $n$ be a
 nonnegative integer.
 Then $\nu_2(n!)=n-s(n)$.
 \end{prop}
 \begin{proof}
 We prove the proposition by induction.
 Clearly $\nu_2(1)=0
 =0-s(0)=1-s(1)$.  Assume that the proposition holds
 for $n\le m-1$.  Let $n=m$.  Suppose first that $m=2l$.
 Then $\nu_2((2l)!)=\nu_2((2l)!!)=
 \nu_2(2^ll!)=l+\nu_2(l!)=l+l-s(l)=m-s(m).$
 Assume now that $m=2l+1$.  Then $\nu_2(m!)=\nu_2((2l)!)=2l-
 s(2l)=m-s(m)$.
\end{proof}

 In what follows we are going to make
 extensive use of the following lemma
 and particularly of its corollary.
 \begin{lemma}\label{coreqfthm1}
 Let $p$ and $q$ be nonnegative integers,
 and assume
 that $n-q= 2p$.  Then
 \begin{equation}\label{forthqn1}
 \nu_2(\th_{q,n})=-p + \sum_{k=1}^p s(n-k)-\sum_{k=1}^p s(k-1).
 \end{equation}
 If $p=0$, we have to interpret empty sums as $0$, as before.
 The equation holds also for $q=0$ if we interpret, as earlier,
 $\th_{0,2p}$ as $1$
 for any $p$.
 \end{lemma}

 \begin{proof}
 Combine
 (\ref{forthqn}) and Proposition~\ref{2adn!}.
 \end{proof}

 Since the $2$-adic valuation in (\ref{forthqn1}) must always be nonnegative,
 we obtain the following corollary.
 \begin{corol}
 For all nonnegative integers $l$ and $p$,
 we have
 \begin{equation}
 \sum_{j=l}^{l+p-1} s(j) \ge p +\sum_{j=0}^{p-1}s(j) \quad \text{if\/ } l\ge
 p,
 \label{basins1}\\
 \end{equation}
 and
 \begin{equation}
 \sum_{j=l}^{l+p-1} s(j)
 \ge l +\sum_{j=0}^{p-1}s(j) \quad \text{if\/ }
 l\le p. \label{basins2}
 \end{equation}
 For $l=p+1,p$, there holds equality,
 \begin{equation}\label{fundeq}
 0= -\sum _{j=0} ^{p-1}s(j)+\sum _{j=p+1} ^{2p}s(j)-p
 =-\sum _{j=0} ^{p-1}s(j)+\sum _{j=p} ^{2p-1}s(j)-p.
 \end{equation}
 In particular,
 \begin{equation} \label{eq:Ungl}
 \sum_{j=l}^{l+p-1} s(j) \ge \sum_{j=0}^{p-1}s(j),
 \end{equation}
 and equality holds if and only if either $p=0$ or $l=0$.
 \end{corol}
 \begin{proof}
 Use (\ref{forthqn1}) with $l=n-p=p+q$ to deduce
 (\ref{basins1}).
 Assume that $0\le l \le p$.  Then by cancelling out
 the common terms in
 (\ref{basins2}) we deduce that
 (\ref{basins2}) follows from
 (\ref{basins1}) with the roles of $l$ and $p$ being interchanged.

 Put $q=0$ in (\ref{forthqn1}) and recall that $\th_{0,n}=1$. This
 implies the second part of (\ref{fundeq}).
 Use the equality $s(2p)=s(p)$ to deduce the first part of
 (\ref{fundeq}).

 The inequality (\ref{eq:Ungl}) follows trivially from
  (\ref{basins1}) and  (\ref{basins2}).
\end{proof}

 \section{Proof of the Falikman--Friedland--Loewy Conjecture}

 In this section, we use the results from the previous section to prove
 the conjecture from \cite{FFL} characterizing the values of $q$ and
 $n$ for which $\delta_{n-q,n}=\th_{q,n}$ is odd. For the sake of
 convenience, we state the result in form of the following theorem.
 The ``if" direction was already shown in \cite{FFL}. Our proof will
 not only show the ``only if" direction, but, in passing, it will also
 provide an independent proof of the ``if" direction.

 \begin{theo} \label{thm:2}
 For positive integers $q$ and $n$, the quantity
 $\th_{q,n}$ is odd if and only if $n\ge q\ge 1$ and
 $$n\equiv \pm q\pmod {2^{\cl{\log_2 2q}}}.$$
 \end{theo}

 \begin{proof}

 For $n-q$ odd, the theorem follows
 from Proposition~\ref{prop:1}. Therefore, for the rest of the
 proof, let $n-q$ be even.

 We repeatedly use
 subsequently the following observation.
 Let $r,t,j$ be nonnegative integers such that $2^t> j$.
 Then $s(r2^t+j)=s(r)+s(j)$.  We divide the proof into the two
 following cases.

 \medskip\noindent
 \textsc{Case 1}.  $n=2n_1$.  It is enough to assume that $q=2q_1$
 and $n_1\ge q_1$.
 Write  $p=n_1-q_1$ and substitute in (\ref{forthqn1}) to
 obtain
 \begin{equation} \label{AE}
 \nu_2(\th_{q,n})=-(n_1-q_1) -\sum _{j=0} ^{n_1-q_1-1}s(j)
 +\sum _{j=n_1+q_1} ^{2n_1-1}s(j).
 \end{equation}
 We now show that the right-hand side of \eqref{AE} is zero if and only if
 $$n_1\equiv \pm q_1\pmod {2^{\cl{\log_2 (2q_1)}}},$$
 which is obviously equivalent to the theorem in this case.

 In \textsc{Case 1} we always use the abbreviation $Q=2^{\cl{\log_2 2q_1}}$.
 Write $n_1=cQ+q_1+d$, where $0\le d<Q$.
 We know that the
 quantity from (\ref{AE}),
 \begin{equation} \label{AK}
 -\sum _{j=0} ^{n_1-q_1-1}s(j)
 +\sum _{j=n_1+q_1} ^{2n_1-1}s(j)-(n_1-q_1)=
 -\sum _{j=0} ^{cQ+d-1}s(j)
 +\sum _{j=cQ+2q_1+d} ^{2cQ+2q_1+2d-1}s(j)-(cQ+d),
 \end{equation}
 is nonnegative.  We have to show that it is zero only if $d=0$ or $d=Q-2q_1$. To do so,
 we distinguish various subcases, depending on the size of $d$.

 \medskip\noindent
 \textsc{Case 1a}: $2q_1+2d\le Q$. In this case, the quantity \eqref{AK} becomes
 \begin{align*}
 -\sum _{r=0} ^{c-1}&\sum _{j=0} ^{Q-1}s(rQ+j)
 -\sum _{j=cQ} ^{cQ+d-1}s(j)
 +\sum _{j=cQ+2q_1+d} ^{(c+1)Q-1}s(j)\\
 &\kern2cm
 +\sum _{r=c+1} ^{2c-1}\sum _{j=0} ^{Q-1}s(rQ+j)
 +\sum _{j=2cQ} ^{2cQ+2q_1+2d-1}s(j)
 -(cQ+d)\\
 &= -Q\sum _{r=0} ^{c-1}s(r)-c\sum _{j=0} ^{Q-1}s(j)
 -ds(c)-\sum _{j=0} ^{d-1}s(j)
 +(Q-2q_1-d)s(c)+\sum _{j=2q_1+d} ^{Q-1}s(j)\\
 &\kern.5cm
 +Q\sum _{r=c+1} ^{2c-1}s(r)+(c-1)\sum _{j=0} ^{Q-1}s(j)
 +(2q_1+2d)s(2c)+\sum _{j=0} ^{2q_1+2d-1}s(j)
 -(cQ+d)\\
 &=
 Q\left(-\sum _{r=0} ^{c-1}s(r)+\sum _{r=c+1} ^{2c}s(r)-c\right)
 -\sum _{j=0} ^{d-1}s(j)
 +\sum _{j=2q_1+d} ^{2q_1+2d-1}s(j) -d.
 \end{align*}
 Using (\ref{fundeq}), we see that the above expression is equal to
 \begin{equation} \label{AL}
 -\sum _{j=0} ^{d-1}s(j)
 +\sum _{j=2q_1+d} ^{2q_1+2d-1}s(j) -d.
 \end{equation}
 By the definitions of $Q$ and $d$, we have $Q/2< 2q_1<
  2q_1+2d-1< Q$. Thus, we have
 \begin{equation} \label{AM}
 -\sum _{j=0} ^{d-1}s(j)
 +\sum _{j=2q_1+d} ^{2q_1+2d-1}s(j)-d
 = -\sum _{j=0} ^{d-1}s(j)
 +\sum _{j=2q_1+d-\frac {Q} {2}} ^{2q_1+2d-\frac {Q} {2}-1}s(j).
 \end{equation}
 It should be noted that, by the definitions of
 $Q$ and $d$, we have $Q/2<2q_1+d$.
 By
 (\ref{eq:Ungl}), the quantity on the right-hand side is
 zero only if the sums on the right-hand side are empty,
 i.e., if $d=0$.

 \medskip\noindent
 \textsc{Case 1b}: $2q_1+d\le Q<2q_1+2d$. In this case, the quantity \eqref{AK} becomes
 \begin{align*}
 -\sum _{r=0} ^{c-1}& \sum _{j=0} ^{Q-1}s(rQ+j)
 -\sum _{j=cQ} ^{cQ+d-1}s(j)
 +\sum _{j=cQ+2q_1+d} ^{(c+1)Q-1}s(j)\\
 &\kern2cm
 +\sum _{r=c+1} ^{2c}\sum _{j=0} ^{Q-1}s(rQ+j)
 +\sum _{j=(2c+1)Q} ^{2cQ+2q_1+2d-1}s(j)
 -(cQ+d)\\
 &= -Q\sum _{r=0} ^{c-1}s(r)-c\sum _{j=0} ^{Q-1}s(j)
 -ds(c)-\sum _{j=0} ^{d-1}s(j)
 +(Q-2q_1-d)s(c)+\sum _{j=2q_1+d} ^{Q-1}s(j)\\
 &\kern.5cm
 +Q\sum _{r=c+1} ^{2c}s(r)+c\sum _{j=0} ^{Q-1}s(j)
 +(2q_1+2d-Q)s(2c)+\sum _{j=Q} ^{2q_1+2d-1}s(j)
 -(cQ+d)\\
 &=
 Q\left(-\sum _{r=0} ^{c-1}s(r)+\sum _{r=c+1} ^{2c}s(r)-c\right)
 -\sum _{j=0} ^{d-1}s(j) +\sum _{j=2q_1+d} ^{2q_1+2d-1}s(j) -d.
 \end{align*}
 Using (\ref{fundeq}) again, we deduce that the
 above expression is equal to
 \begin{equation} \label{AN}
 -\sum _{j=0} ^{d-1}s(j)
 +\sum _{j=2q_1+d} ^{2q_1+2d-1}s(j) -d.
 \end{equation}
 By the definitions of $Q$ and $d$, we have $Q/2< 2q_1<
 2q_1+d\le Q\le 2q_1+2d-1<2Q$. Thus, the above quantity can be further
 modified to
 \begin{align*}
 -\sum _{j=0} ^{d-1}s(j)
 +\sum _{j=2q_1+d} ^{Q-1}s(j)
 +&\sum _{j=Q} ^{2q_1+2d-1}s(j) -d\\
 &=-\sum _{j=0} ^{d-1}s(j)
 +\sum _{j=2q_1+d-\frac {Q} {2}} ^{\frac {Q} {2}-1}s(j)
 +\sum _{j=0} ^{2q_1+2d-Q-1}s(j)
 \\
 &=-\sum _{j=2q_1+2d-Q} ^{d-1}s(j)
 +\sum _{j=2q_1+d-\frac {Q} {2}} ^{\frac {Q} {2}-1}s(j).
 \end{align*}
 From $Q/2<2q_1$ and $2q_1+d\le Q$, we infer that $d<Q/2$.
 Now we use identity \eqref{AI} with $2^e=\frac{Q}{2}$ for all the digit sums
 in the last expression. This leads to the expression
 $$\sum _{j=\frac {Q} {2}-d} ^{\frac {3} {2}Q-2q_1-2d-1}s(j)
 -\sum _{j=0} ^{Q-2q_1-d-1}s(j). $$
 We recall that $Q-2q_1-d\ge 0$.
 Apply
 (\ref{eq:Ungl}) again to conclude
 that the last expression, and hence, $\nu_2(\th_{q,n})$, is
 zero only if the sums in the last line are empty,
 that is, if $d=Q-2q_1$.

 \medskip\noindent
 \textsc{Case 1c}: $Q<2q_1+d$ and $2q_1+2d\le 2Q$.
 In this case, the quantity \eqref{AK} becomes
 \begin{align} \notag
 -\sum _{r=0} ^{c-1}&\sum _{j=0} ^{Q-1}s(rQ+j)
 -\sum _{j=cQ} ^{cQ+d-1}s(j)
 +\sum _{j=cQ+2q_1+d} ^{(c+2)Q-1}s(j)\\
 \notag
 &\kern2cm
 +\sum _{r=c+2} ^{2c}\sum _{j=0} ^{Q-1}s(rQ+j)
 +\sum _{j=(2c+1)Q} ^{2cQ+2q_1+2d-1}s(j)
 -(cQ+d)\\
 \notag
 &= -Q\sum _{r=0} ^{c-1}s(r)-c\sum _{j=0} ^{Q-1}s(j)
 -ds(c)-\sum _{j=0} ^{d-1}s(j)
 +(2Q-2q_1-d)s(c+1)\\
 \notag
 &\kern2cm
  +\sum _{j=2q_1+d-Q} ^{Q-1}s(j)
 +Q\sum _{r=c+2} ^{2c}s(r)+(c-1)\sum _{j=0} ^{Q-1}s(j)\\
 \notag
 &\kern2cm
 +(2q_1+2d-Q)s(2c+1)+\sum _{j=0} ^{2q_1+2d-Q-1}s(j) -(cQ+d)\\
 \notag
 &= Q\left(-\sum _{r=0} ^{c-1}s(r)+\sum _{r=c+1} ^{2c}s(r)-c\right)
 -\sum _{j=0} ^{d-1}s(j)
 +\sum _{j=2q_1+d-Q} ^{2q_1+2d-Q-1}s(j)\\
 \notag
 &\kern2cm
 +(2q_1+d-Q)(s(c)-s(c+1)+1)\\
 &= -\sum _{j=0} ^{d-1}s(j)+\sum _{j=2q_1+d-Q} ^{2q_1+2d-Q-1}s(j)
 +(2q_1+d-Q)(s(c)-s(c+1)+1).
 \label{ANa}
 \end{align}
 In the second step we used the equality $s(2c+1)=s(c)+1$, and
 in the last step we used again (\ref{fundeq}).
 Since $Q<2q_1+d$, we have $2q_1+d-Q>0$. As
 $s(c)-s(c+1)+1\ge0$ for any $c\ge 0$, the third term in the last line is
 nonnegative.   As $d>0$,
 the inequality (\ref{eq:Ungl}) implies that the sum of the first two terms is
 strictly positive.

 \medskip\noindent
 \textsc{Case 1d}: $Q<2q_1+d$ and $2q_1+2d> 2Q$.
 In this case, the quantity \eqref{AK} becomes
 \begin{align*}
 -\sum _{r=0} ^{c-1}&\sum _{j=0} ^{Q-1}s(rQ+j)
 -\sum _{j=cQ} ^{cQ+d-1}s(j)
 +\sum _{j=cQ+2q_1+d} ^{(c+2)Q-1}s(j)\\
 &\kern2cm
 +\sum _{r=c+2} ^{2c+1}\sum _{j=0} ^{Q-1}s(rQ+j)
 +\sum _{j=(2c+2)Q} ^{2cQ+2q_1+2d-1}s(j)
 -(cQ+d)\\
 &=
 -Q\sum _{r=0} ^{c-1}s(r)-c\sum _{j=0} ^{Q-1}s(j)
 -ds(c)-\sum _{j=0} ^{d-1}s(j)
 +(2Q-2q_1-d)s(c+1)\\
 &\kern2cm
 +\sum _{j=2q_1+d-Q} ^{Q-1}s(j)
 +Q\sum _{r=c+2} ^{2c+1}s(r)+c\sum _{j=0} ^{Q-1}s(j)\\
 &\kern2cm
 +(2q_1+2d-2Q)s(2c+2)+\sum _{j=0} ^{2q_1+2d-2Q-1}s(j)
 -(cQ+d).
 \end{align*}
 We now do the following substitutions.  First, $s(2(c+1))=s(c+1)$.
 Second, in the sum
 over $j=0,\ldots,2q_1+2d-2Q-1$
 (where $2q_1+2d-2Q-1\le
 Q+2(Q-1)-2Q-1=Q-1$), we let $s(j)=s(Q+j)-1$.
 Third,
 $$\sum _{r=c+2} ^{2c+1}s(r)=s(2c+1)-s(c+1)+\sum _{r=c+1}
 ^{2c}s(r)=s(c)+1-s(c+1) +\sum _{r=c+1} ^{2c}s(r).$$
 Hence, the above expression is equal to
 \begin{align}
 \notag
 &Q\left(-\sum _{r=0} ^{c-1}s(r)+\sum _{r=c+1} ^{2c}s(r)-c\right)
 -\sum _{j=0} ^{d-1}s(j)
 +\sum _{j=2q_1+d-Q} ^{Q-1}s(j)\\
 \notag
 &\kern2cm
 +\sum _{j=Q} ^{2q_1+2d-Q-1}s(j)-(2q_1+2d-2Q)+
 (Q-d)(s(c)-s(c+1)+1)\\
 &=-\sum _{j=0} ^{d-1}s(j)+\sum _{j=2q_1+d-Q} ^{2q_1+2d-Q-1}s(j)
 -(2q_1+2d-2Q)+(Q-d)(s(c)-s(c+1)+1),
 \label{ANb}
 \end{align}
 where we used again (\ref{fundeq}).
 Since $Q>d$ and $s(c)-s(c+1)+1\ge 0$ for any $c\ge0$, the fourth term in
 the last line is nonnegative.
 We have $\frac{Q}{2}<d<Q$ and $d\ge 2q_1+d-Q$.
 Thus, we may apply (\ref{basins2})
 to deduce that the sum of the first two terms is at least
 $2q_1+d-Q$.  Hence the expression (\ref{ANb}) is not less than
 $Q-d>0$.

 The proof of \textsc{Case~1} is completed.

 \medskip\noindent
 \textsc{Case 2}.  Let $n=2n_1+1$, $q=2q_1+1$, where $n_1\ge q_1\ge 0$.
 Write $p=n_1-q_1$ and substitute in (\ref{forthqn1}) to obtain
 \begin{equation} \label{AF}
 \nu_2(\th_{q,n})=-(n_1-q_1)
 -\sum _{j=0} ^{n_1-q_1-1}s(j)
 +\sum _{j=n_1+q_1+1} ^{2n_1}s(j).
 \end{equation}
 We now show that the right-hand side of \eqref{AF} is zero if and only if
 $$n_1\equiv  q_1\pmod {2^{\cl{\log_2 (2q_1+1)}}}$$
 or
 $$n_1\equiv  -q_1-1\pmod {2^{\cl{\log_2 (2q_1+1)}}},$$
 which is obviously equivalent to the theorem in this case.

 In \textsc{Case~2} we always use the abbreviation
 $Q=2^{\cl{\log_2 (2q_1+1)}}$.
 Write $n_1=cQ+q_1+d$, where $0\le d<Q$.
 We have to show that
 the expression (\ref{AF}), that is,
 \begin{equation}
 -\sum _{j=0} ^{n_1-q_1-1}s(j)
 +\sum _{j=n_1+q_1+1} ^{2n_1}s(j)-(n_1-q_1)
 = -\sum _{j=0} ^{cQ+d-1}s(j)
 +\sum _{j=cQ+2q_1+d+1} ^{2cQ+2q_1+2d}s(j)-(cQ+d),
 \label{AO}
 \end{equation}
 is zero only if $d=0$ or $d=Q-2q_1-1$. To do so,
 we distinguish again various cases, depending on the size of $d$.

 \medskip\noindent
 \textsc{Case 2a}: $2q_1+2d< Q$.
 In this case, the quantity \eqref{AO} becomes
 \begin{align*}
 -{}&\sum _{r=0} ^{c-1}\sum _{j=0} ^{Q-1}s(rQ+j)
 -\sum _{j=cQ} ^{cQ+d-1}s(j)
 +\sum _{j=cQ+2q_1+d+1} ^{(c+1)Q-1}s(j)\\
 &\kern2cm
 +\sum _{r=c+1} ^{2c-1}\sum _{j=0} ^{Q-1}s(rQ+j)
 +\sum _{j=2cQ} ^{2cQ+2q_1+2d}s(j)
 -(cQ+d)\\
 &=
 -Q\sum _{r=0} ^{c-1}s(r)-c\sum _{j=0} ^{Q-1}s(j)
 -ds(c)-\sum _{j=0} ^{d-1}s(j)
 +(Q-2q_1-d-1)s(c)+\sum _{j=2q_1+d+1} ^{Q-1}s(j)\\
 &\kern.5cm
 +Q\sum _{r=c+1} ^{2c-1}s(r)+(c-1)\sum _{j=0} ^{Q-1}s(j)
 +(2q_1+2d+1)s(2c)+\sum _{j=0} ^{2q_1+2d}s(j)
 -(cQ+d)\\
 &=
 Q\left(-\sum _{r=0} ^{c-1}s(r)+\sum _{r=c+1} ^{2c}s(r)-c\right)
 -\sum _{j=0} ^{d-1}s(j)
 +\sum _{j=2q_1+d+1} ^{2q_1+2d}s(j) -d\\
 &=
 -\sum _{j=0} ^{d-1}s(j)
 +\sum _{j=2q_1+d+1} ^{2q_1+2d}s(j) -d,
 \end{align*}
 where in the last line we used again (\ref{fundeq}).
 The reader should compare this expression with the one in \eqref{AL}.
 Indeed, the remaining arguments are completely analogous to those
 after \eqref{AL} in \textsc{Case~1a} of the current proof, which are
 therefore left to the reader.

 \medskip\noindent
 \textsc{Case 2b}: $2q_1+d< Q\le 2q_1+2d$.
 In this case, the quantity \eqref{AO} becomes
 \begin{align*}
 -{}&\sum _{r=0} ^{c-1}\sum _{j=0} ^{Q-1}s(rQ+j)
 -\sum _{j=cQ} ^{cQ+d-1}s(j)
 +\sum _{j=cQ+2q_1+d+1} ^{(c+1)Q-1}s(j)\\
 &\kern2cm
 +\sum _{r=c+1} ^{2c}\sum _{j=0} ^{Q-1}s(rQ+j)
 +\sum _{j=(2c+1)Q} ^{2cQ+2q_1+2d}s(j)
 -(cQ+d)\\
 &=
 -Q\sum _{r=0} ^{c-1}s(r)-c\sum _{j=0} ^{Q-1}s(j)
 -ds(c)-\sum _{j=0} ^{d-1}s(j)
 +(Q-2q_1-d-1)s(c)+\sum _{j=2q_1+d+1} ^{Q-1}s(j)\\
 &\kern.5cm
 +Q\sum _{r=c+1} ^{2c}s(r)+c\sum _{j=0} ^{Q-1}s(j)
 +(2q_1+2d+1-Q)s(2c)+\sum _{j=Q} ^{2q_1+2d}s(j)
 -(cQ+d)\\
 &=
 Q\left(-\sum _{r=0} ^{c-1}s(r)+\sum _{r=c+1} ^{2c}s(r)-c\right)
 -\sum _{j=0} ^{d-1}s(j)
 +\sum _{j=2q_1+d+1} ^{2q_1+2d}s(j) -d\\
 &= -\sum _{j=0} ^{d-1}s(j)
 +\sum _{j=2q_1+d+1} ^{2q_1+2d}s(j) -d,
 \end{align*}
 where in the last line we used again (\ref{fundeq}).
 The reader should compare this expression with the one in \eqref{AN}.
 Indeed, the remaining arguments are completely analogous to those
 after \eqref{AN} in \textsc{Case~1b} of the current proof, which are
 therefore left to the reader.

 \medskip\noindent
 \textsc{Case 2c}. $Q\le 2q_1+d$ and $2q_1+2d<2Q$.
 In this case, the quantity \eqref{AO} becomes
 \begin{align*}
 -\sum _{r=0} ^{c-1}&\sum _{j=0} ^{Q-1}s(rQ+j)
 -\sum _{j=cQ} ^{cQ+d-1}s(j)
 +\sum _{j=cQ+2q_1+d+1} ^{(c+2)Q-1}s(j)\\
 &\kern2cm
 +\sum _{r=c+2} ^{2c}\sum _{j=0} ^{Q-1}s(rQ+j)
 +\sum _{j=(2c+1)Q} ^{2cQ+2q_1+2d}s(j)
 -(cQ+d)\\
 &=
 -Q\sum _{r=0} ^{c-1}s(r)-c\sum _{j=0} ^{Q-1}s(j)
 -ds(c)-\sum _{j=0} ^{d-1}s(j)
 +(2Q-2q_1-d-1)s(c+1)\\
 &\kern2cm
 +\sum _{j=2q_1+d+1-Q} ^{Q-1}s(j)
 +Q\sum _{r=c+2} ^{2c}s(r)+(c-1)\sum _{j=0} ^{Q-1}s(j)\\
 &\kern2cm
 +(2q_1+2d+1-Q)s(2c+1)+\sum _{j=0} ^{2q_1+2d-Q}s(j)
 -(cQ+d)\\
 &= Q\left(-\sum _{r=0} ^{c-1}s(r)+\sum _{r=c+1} ^{2c}s(r)-c\right)
 -\sum _{j=0} ^{d-1}s(j) +\sum _{j=2q_1+d+1-Q} ^{2q_1+2d-Q}s(j)\\
 &\kern2cm
 +(2q_1+d+1-Q)(s(c)-s(c+1)+1)\\
 &= -\sum _{j=0} ^{d-1}s(j)
 +\sum _{j=2q_1+d+1-Q} ^{2q_1+2d-Q}s(j)
 +(2q_1+d+1-Q)(s(c)-s(c+1)+1),
 \end{align*}
 where in the last step we used again (\ref{fundeq}).
 Since $Q<2q_1+d+1$, we have $2q_1+d+1-Q>0$. In particular, since also
 $s(c)-s(c+1)+1\ge0$ for any $c$, the third term in the last line is
 nonnegative, and, because $d>0$,
 the inequality (\ref{eq:Ungl}) says that the sum of the first two terms is
 strictly positive.

 \medskip\noindent
 \textsc{Case 2d}. $Q\le 2q_1+d$ and $2q_1+2d\ge 2Q$.
 In this case, the quantity \eqref{AO} becomes
 \begin{align*}
 -&\sum _{r=0} ^{c-1}\sum _{j=0} ^{Q-1}s(rQ+j)
 -\sum _{j=cQ} ^{cQ+d-1}s(j)
 +\sum _{j=cQ+2q_1+d+1} ^{(c+2)Q-1}s(j)\\
 &\kern2cm
 +\sum _{r=c+2} ^{2c+1}\sum _{j=0} ^{Q-1}s(rQ+j)
 +\sum _{j=(2c+2)Q} ^{2cQ+2q_1+2d}s(j)
 -(cQ+d)\\
 &=
 -Q\sum _{r=0} ^{c-1}s(r)-c\sum _{j=0} ^{Q-1}s(j)
 -ds(c)-\sum _{j=0} ^{d-1}s(j)
 +(2Q-2q_1-d-1)s(c+1)\\
 &\kern2cm
 +\sum _{j=2q_1+d+1-Q} ^{Q-1}s(j)
 +Q\sum _{r=c+2} ^{2c+1}s(r)+c\sum _{j=0} ^{Q-1}s(j)\\
 &\kern2cm
 +(2q_1+2d+1-2Q)s(2c+2)+\sum _{j=0} ^{2q_1+2d-2Q}s(j)
 -(cQ+d)\\
 &=
 Q\left(-\sum _{r=0} ^{c-1}s(r)+\sum _{r=c+1} ^{2c}s(r)-c\right)
 -\sum _{j=0} ^{d-1}s(j)
 +\sum _{j=2q_1+d+1-Q} ^{Q-1}s(j)\\
 &\kern2cm
 +\sum _{j=Q} ^{2q_1+2d-Q}s(j)
 -(2q_1+2d-2Q+1)
 +(Q-d)(s(c)-s(c+1)+1)\\
 &= -\sum _{j=0} ^{d-1}s(j)
 +\sum _{j=2q_1+d+1-Q} ^{2q_1+2d-Q}s(j)
 -(2q_1+2d-2Q+1)
 +(Q-d)(s(c)-s(c+1)+1),
 \end{align*}
 where in the last step we used again (\ref{fundeq}).
 We have $d<Q$. In particular, since also
 $s(c)-s(c+1)+1\ge0$ for any $c$, the fourth term in the last line is
 nonnegative. Moreover, because $2q_1+d+1-Q\le d$,
 the inequality (\ref{basins2})
 says that the sum of the first two terms is
 at least $2q_1+d+1-Q=(2q_1+2d-2Q+1)+(Q-d)$. Thus, the sum of the first
 three terms in the last line is strictly positive. Combining both findings,
 we infer that the whole last line is strictly positive
 (in fact, at least $Q-d$).

 This finishes the proof of \textsc{Case~2} and, thus, of the theorem.
 \end{proof}

 \section{Additional results on the 2-adic
 valuation of $\delta_{n-q,n}=\th_{q,n}$}

 In this section, we examine the $2$-adic behavior of
 $\th_{q,n}=\delta_{n-q,n}$ in more detail. Keeping in mind
 Proposition~\ref{prop:1}, we concentrate throughout this section on
 the case that $n$ and $q$ have the same parity.

 If we fix $q$ and let $n=q+2i$, $i=0,1,2,\dots$, vary, then we know
 that whenever (\ref{spc1}) is satisfied, the $2$-adic valuation of
 $\th_{q,n}=\th_{q,q+2i}$ will be zero. However, what happens in between?
 By looking at some random values of $q$, one might get the impression
 that, between two successive occurrences of zero, the $2$-adic valuations
 $\nu_2(\th_{q,q+2i})$ are unimodal, that is, they
 first grow (weakly) monotone until they reach their maximum value
 half-way, and then they drop (weakly) monotone until they come back to zero in the
 end. Moreover, one is led to guess that the $2$-adic valuations are
 symmetric around the place where the maximum is attained. As it turns
 out (see Example~\ref{Ex1} below), the unimodality conjecture is not
 true, while the symmetry conjecture is indeed true. However, in some
 sense, unimodality is ``almost" true. As we demonstrate in
 Theorems~\ref{thm:8}, \ref{thm:9}, \ref{thm:10}, \ref{thm:11} below,
 between two successive values of zero, the $2$-adic valuations
 $\nu_2(\th_{q,q+2i})$ stay above a linear function of slope $1$ which
 is tight at the opening zero, and at the same time they stay below
 another linear function of slope 1 which is tight at the maximum
 (see (\ref{AP}) and (\ref{AQ}) in Theorem~\ref{thm:8} and the analogous
 inequalities in the subsequent theorems),
 until they reach the maximum value, which is attained exactly
 half-way, and the $2$-adic valuations beyond are the symmetric
 images of those before. In the theorems, we determine in addition
 the maximal value for each of these intervals.

 The results which are found on the way to prove these theorems allow
 one also to address the following question: characterize all values of $q$
 and $n$ for which the $2$-adic valuation of $\th_{q,n}$ has a certain
 fixed (small) value.
 Clearly, Theorem~\ref{thm:2} does this if we fix this value to
 $0$. In Corollaries~\ref{cor:8}, \ref{cor:9}, \ref{cor:10},
 \ref{cor:11} we work out the analogous
 characterization if we fix this value to
 $1$. We could move on to $2$, $3$, etc., but at the cost of
 a considerable increase of complication the further we go.

 \begin{example} \label{Ex1}
 The sequence $(\nu_2(\th_{q,q+2i}))_{i\ge0}$ is not unimodal in the
 intervals discussed in Theorems~\ref{thm:8}--\ref{thm:11}
 (although the theorems show that
 it comes very close).
 Here we give two counter-examples.

 Consider $q=39$. In this case, the first
 200 terms of the sequence  $(\nu_2(\th_{q,q+2i}))_{i\ge0}$ are:
 \begin{align*}
 0, &1, 3, 5, 7, 8, 8, 8, 8, 9, 12, 15, 18, 18, 15, 12, 9, 8, 8, 8, 8,
     7, 5, 3, 1, 0, 3, 6, 9, 10, 9, 8, 7, 9, 12, 15,\\
 & 18, 20, 21, 22, 23, 25,
     29, 33, 37, 37, 33, 29, 25, 23, 22, 21, 20, 18, 15, 12, 9, 7, 8, 9, 10,
     9, 6, 3,\\
 & 0, 1, 3, 5, 7, 8, 8, 8, 8, 9, 12, 15, 18, 18, 15, 12, 9, 8, 8,
     8, 8, 7, 5, 3, 1, 0, 4, 8, 12, 14, 14, 14, 14,\\
 & 17, 21, 25, 29, 32, 34,
     36, 38, 41, 46, 51, 56, 56, 51, 46, 41, 38, 36, 34, 32, 29, 25, 21, 17,
     14, 14, \\
 &14, 14, 12, 8, 4, 0, 1, 3, 5, 7, 8, 8, 8, 8, 9, 12, 15, 18, 18,
     15, 12, 9, 8, 8, 8, 8, 7, 5, 3, 1, 0, 3, 6, 9, 10,\\
 & 9, 8, 7, 9, 12, 15,
     18, 20, 21, 22, 23, 25, 29, 33, 37, 37, 33, 29, 25, 23, 22, 21, 20, 18,
     15, 12, 9, 7, \\
 & 8,9, 10, 9, 6, 3, 0, 1, 3, 5, 7, 8, 8, 8, 8.
 \end{align*}
 As predicted by Theorem~\ref{thm:2}, zeros occur for
 all multiples of $\frac {1} {2}2^{\cl{\log_22\cdot 39}}=64$ and for
 multiples of $64$ reduced by $39$. However, while the first
 interval (the interval comprising the first 25 values) is unimodal,
 the second one (comprising the next 39 values) is not, as can be seen
 from the subsequence $0, 3, 6, 9, 10, 9, 8, 7, 9, 12, 15,\dots, 37,\dots$,
 which rises first to $10$, then drops
 to $7$, to rise again beyond $10$ until it reaches the maximum of
 $37$ in this interval.

 Consider now $q=46$. In this case, the first 200 terms are:
 \begin{align*} 0,& 4, 2, 5, 6, 10, 10, 13, 14, 19, 14, 13, 10, 10, 6, 5, 2, 4, 0, 3,
     4, 8, 8, 11, 12, 17, 14, 15, 14, 16, 14,\\
 & 15, 14, 19, 18, 22, 24, 29, 30,
     34, 36, 42, 36, 34, 30, 29, 24, 22, 18, 19, 14, 15, 14, 16, 14,
 15, 14,\\
 &
     17, 12, 11, 8, 8, 4, 3, 0, 4, 2, 5, 6, 10, 10, 13, 14, 19, 14, 13, 10,
     10, 6, 5, 2, 4, 0, 4, 6, 11, 12, 16,  \\
 &18, 24, 22, 24, 24, 27, 26, 28, 28,
     34, 34, 39, 42, 48, 50, 55, 58, 65, 58, 55, 50, 48, 42, 39, 34, 34, \\
 &28, 28, 26, 27, 24, 24, 22, 24, 18, 16, 12, 11, 6, 4, 0, 4, 2, 5, 6, 10, 10,
     13, 14, 19, 14, 13, 10, 10, 6,   \\
 &5, 2, 4, 0, 3, 4, 8, 8, 11, 12, 17, 14,
     15, 14, 16, 14, 15, 14, 19, 18, 22, 24, 29, 30, 34, 36, 42, 36, 34,
    \\
 & 30, 29, 24, 22, 18, 19, 14, 15, 14, 16, 14, 15, 14, 17, 12, 11, 8, 8, 4, 3,
     0, 4, 2, 5, 6, 10, 10, 13, 14.
 \end{align*}
%Christian: correction
 This is an example where {\it both\/} types of intervals are not
 unimodular. The first interval (comprising the first $64-46=18$ values)
 begins $ 0, 4, 2, 5, 6, 10, 10, 13, 14, 19, \dots$,
 that is, rises to $4$, drops 
 to $2$, before rising again up to the maximum of $19$ in this interval.
 The second one (comprising the next $46$ values) starts by
 $0, 3,  4, 8, 8, 11, 12, 17, 14, 15,\dots,42, \dots$, that is, rises to
 $17$, drops to $14$, before rising again up to the maximum of $42$
 in this interval. Similar remarks apply to the subsequent intervals.
 \end{example}

 Nevertheless, between two successive zeros, although the $2$-adic valuations
 $\nu_2(\th_{q,q+2i})$ are not unimodal in general, the $2$-adic valuations still seem
 to exhibit an overall increase until a maximum halfway and then
 a decrease which is the symmetric image of the increasing values.
 In the theorems below, we quantify this statement.
 We split our results into four separate theorems.
 On the one hand, we have to distinguish between the two possible
 parities of $q$, and, on the other hand, for each integer $q$, there
 are two types of intervals to be considered.

 \begin{theo} \label{thm:8}  Let $q$ be a fixed positive even integer,
 $Q=2^{\cl{\log_2 q}}$, and let $c$ be any nonnegative integer.
 Then the values of the sequence
 $(\nu_2(\th_{q,q+2i}))_{i\ge0}$ for
 $i\in [cQ,(c+1)Q-q]$ are symmetric around its center $i_1=cQ+(Q-q)/2$.
 The values at the extreme points of the interval
 $i=cQ$ and $i=(c+1)Q-q$ are $0$, and
 the value at the center is
 $$\nu_2(\th_{q,q+2i_1})=
 \frac {Q-q} {2}\log_2\frac {Q} {2}
 -2\sum _{j=0} ^{\frac {Q-q-2} {2}}s(j),
 $$
 independent of $c$.  Furthermore,
 for $i\in [cQ,i_1]$ the inequalities
 \begin{equation} \label{AP}
 \nu_2(\th_{q,q+2i})\ge i-cQ
 \end{equation}
 and
 \begin{equation} \label{AQ}
 \nu_2(\th_{q,q+2i})\le \nu_2(\th_{q,q+2i_1})-\(i_1-i\)
 \end{equation}
 hold.
 \end{theo}

 \begin{proof}
 In this proof, we use again the notation $q_1=q/2$.
 To show the symmetry  write
 $ i=cQ+d$ with $0\le d\le Q-2q_1$. In particular, we have $2q_1+d\le
 Q$, and, hence (cf.\ \eqref{AL} and \eqref{AN}), the 2-adic
 valuation of
 $\th_{q,q+2i}$ is given by
 \begin{equation} \label{AR}
 \nu_2(\th_{q,q+2i})=\nu_2(\th_{q,q+2(cQ+d)})
 =-\sum _{j=0} ^{d-1}s(j)
 +\sum _{j=2q_1+d} ^{2q_1+2d-1}s(j) -d.
 \end{equation}
 On the other hand, for the same reason, the 2-adic
 valuation of
 $\th_{q,q+2(cQ+Q-q-d)}$ is given by
 \begin{align*}
 \nu_2(&\th_{q,q+2(cQ+Q-2q_1-d)})=-\sum _{j=0} ^{Q-2q_1-d-1}s(j)
 +\sum _{j=2q_1+(Q-2q_1-d)} ^{2q_1+2(Q-2q_1-d)-1}s(j)
 -(Q-2q_1-d)\\
 &=-\sum _{j=0} ^{Q-2q_1-d-1}s(j)
 +\sum _{j=Q-d} ^{2Q-2q_1-2d-1}s(j)
 -(Q-2q_1-d)\\
 &=-\sum _{j=0} ^{Q-1}s(j)+\sum _{j=Q-2q_1-d} ^{Q-1}s(j)\\
 &\kern2cm
 +\sum _{j=0} ^{2Q-1}s(j)
 -\sum _{j=0} ^{Q-d-1}s(j)
 -\sum _{j=2Q-2q_1-2d} ^{2Q-1}s(j)
 -(Q-2q_1-d)\\
 &=\sum _{j=Q-2q_1-d} ^{Q-1}s(j)
 +\sum _{j=Q} ^{2Q-1}s(j)
 -\sum _{j=0} ^{Q-d-1}s(j)
 -\sum _{j=2Q-2q_1-2d} ^{2Q-1}s(j)
 -(Q-2q_1-d)\\
 &=\sum _{j=Q-2q_1-d} ^{Q-1}s(j)
 +\sum _{j=Q-d} ^{Q-1}s(j)
 -\sum _{j=2Q-2q_1-2d} ^{2Q-1}s(j)
 +(2q_1+d).
 \end{align*}

 Now we apply the reflection identity \eqref{AI} to all the sum of digit
 functions. Thus, we obtain
 \begin{align*}
 \nu_2(\th_{q,q+2(cQ+Q-2q_1-d)})
 &=(2q_1+2d)\log_2Q-\sum _{j=0} ^{2q_1+d-1}s(j)
 -\sum _{j=0} ^{d-1}s(j)\\
 &\kern2cm
 -(2q_1+2d)\log_2(2Q)+\sum _{j=0} ^{2q_1+2d-1}s(j)
 +(2q_1+d)\\
 &=
 -\sum _{j=0} ^{d-1}s(j)
 +\sum _{j=2q_1+d} ^{2q_1+2d-1}s(j)
 -d\\
 &=\nu_2(\th_{q,q+2(cQ+d)}),
 \end{align*}
 proving the symmetry.
 That the values of the extreme points of the interval $i=cQ$ and
 $i=(c+1)Q-q$ are $0$, was already shown in Theorem~\ref{thm:2}, \textsc{Case~1}.

 Next we determine the 2-adic value of $\th_{q,q+2i}$ at the
 center $i_1=cQ+(Q-q)/2$ of the interval. By \eqref{AR}, we have
 \begin{align*}
 \nu_2(\th_{q,q+2i_1})&=
 \nu_2(\th_{q,q+2(cQ+Q/2-q_1)})\\
 &=
 -\sum _{j=0} ^{\frac {Q} {2}-q_1-1}s(j)
 +\sum _{j=\frac {Q} {2}+q_1} ^{Q-1}s(j)
 -\left(\frac {Q} {2}-q_1\right)\\
 &=
 -\sum _{j=0} ^{\frac {Q} {2}-q_1-1}s(j)
 +\sum _{j=q_1} ^{\frac {Q} {2}-1}s(j)\\
 &=
 \(\frac {Q} {2}-q_1\)\log_2\frac {Q} {2}
 -2\sum _{j=0} ^{\frac {Q} {2}-q_1-1}s(j).
 \end{align*}

 To prove the inequality \eqref{AP}, we write again $i=cQ+d$,
 with $d\le Q/2-q_1$. Under this condition, we
 showed in the proof of Theorem~\ref{thm:2}, \textsc{Case~1a} (see \eqref{AM}),
 that $\nu_2(\th_{q,q+2i})$ can be written in the form
 \begin{equation} \label{AP1}
 \nu_2(\th_{q,q+2i})=\nu_2(\th_{q,q+2(cQ+d)})=-\sum _{j=0} ^{d-1}s(j)
 +\sum _{j=2q_1+d-\frac {Q} {2}} ^{2q_1+2d-\frac {Q} {2}-1}s(j).
 \end{equation}
 Since $d\le 2q_1+d-Q$,
 the inequality (\ref{basins1}) implies that this expression is at
 least $d=i-cQ$.

 To prove inequality \eqref{AQ}, we compute the difference of the 2-adic
 valuations of $\th_{q,q+2i_1}$ and $\th_{q,q+2i}$ for
 $i\in[cQ,i_1]$: let again $i=cQ+d$,
 $0\le d\le (Q-q)/2$. Then, using \eqref{AR}, we have
 \begin{align*}
 \nu_2(\th_{q,q+2i_1})-\nu_2(\th_{q,q+2i})&=
 \nu_2(\th_{q,q+2(cQ+Q/2-q_1)})-\nu_2(\th_{q,q+2(cQ+d)})\\
 &=-\sum _{j=0} ^{\frac {Q} {2}-q_1-1}s(j)
 +\sum _{j=\frac {Q} {2}+q_1} ^{Q-1}s(j)
 -\left(\frac {Q} {2}-q_1\right)
\\
&\kern2cm
 +\sum _{j=0} ^{d-1}s(j)
 -\sum _{j=2q_1+d} ^{2q_1+2d-1}s(j)
 +d\\
 &=-\sum _{j=d} ^{\frac {Q} {2}-q_1-1}s(j)
 -\sum _{j=2q_1+d} ^{\frac {Q} {2}+q_1-1}s(j)
 +\sum _{j=2q_1+2d} ^{Q-1}s(j)
 -\left(\frac {Q} {2}-q_1-d\right).
 \end{align*}
 Now we apply the reflection formula \eqref{AI} to all the sum of digits
 functions. This yields
 \begin{align*}
 \nu_2(\th_{q,q+2i_1})-\nu_2(&\th_{q,q+2i})
 =\sum _{j=q_1} ^{\frac {Q} {2}-d-1}s(j)
 +\sum _{j=\frac {Q} {2}-q_1} ^{Q-2q_1-d-1}s(j)
 -\sum _{j=0} ^{Q-2q_1-2d-1}s(j)\\
 &=\sum _{j=q_1} ^{\frac {Q} {2}-d-1}s(j)
 +\sum _{j=\frac {Q} {2}-q_1} ^{Q-2q_1-d-1}s(j)
 -\sum _{j=0} ^{\frac {Q} {2}-q_1-d-1}(s(2j)+s(2j+1))\\
 &=\sum _{j=q_1} ^{\frac {Q} {2}-d-1}s(j)
 +\sum _{j=\frac {Q} {2}-q_1} ^{Q-2q_1-d-1}s(j)
 -2\sum _{j=0} ^{\frac {Q} {2}-q_1-d-1}s(j)
 -\(\frac {Q} {2}-q_1-d\).
 \end{align*}
 Since $q_1\ge \frac {Q} {2}-q_1\ge \frac {Q} {2}-q_1-d$,
 we may apply
 (\ref{basins1}) twice to obtain finally
 $$
 \nu_2(\th_{q,q+2i_1})-\nu_2(\th_{q,q+2i})
 \ge\(\frac {Q} {2}-q_1-d\)=i_1-i,
 $$
 as desired.
\end{proof}

 \begin{corol} \label{cor:8}
 For a fixed even $q$, the values of the sequence
 $(\nu_2(\th_{q,q+2i}))_{i\ge0}$ are never $1$ for
 $i\in [cQ,(c+1)Q-q]$.
 \end{corol}
 \begin{proof}  By the symmetry of the values around the center of the
 interval, and by the inequality \eqref{AP}, the only possible
 candidates are $i=cQ+1$ and $i=(c+1)Q-q-1$. Formula~\eqref{AP1} for
 $\nu_2(\th_{q,q+2i})$ implies that in that case we must have
 $$s\(q+1-\frac {Q} {2}\)=1.$$
 This is only possible if $q+1-\frac {Q} {2}$ is a power of $2$, and
 this, in its turn, implies that $q$ must be odd, a contradiction.
\end{proof}

 \begin{theo} \label{thm:9}
 Let $q$ be a fixed positive even integer,
 $Q=2^{\cl{\log_2 q}}$, and let $c$ be any nonnegative integer.
 Then the values of the sequence
 $(\nu_2(\th_{q,q+2i}))_{i\ge0}$ for
 $i\in [(c+1)Q-q,(c+1)Q]$ are symmetric around its center $i_2=(c+1)Q-q/2$.
 The values at the extreme points of the interval $i=(c+1)Q-q$ and
 $i=(c+1)Q$ are $0$, and the value at the center is
 $$\nu_2(\th_{q,q+2i_2})=
 \frac {q} {2}\log_2Q
 -2\sum _{j=0} ^{\frac {q-2} {2}}s(j)+\frac {q} {2}(s(c)-s(c+1)+1).
 $$
 Furthermore,
 for $i\in [i_2,(c+1)Q]$ the inequalities
 \begin{equation} \label{AT}
 \nu_2(\th_{q,q+2i})\ge (c+1)Q-i
 \end{equation}
 and
 \begin{equation} \label{AS}
 \nu_2(\th_{q,q+2i})\le \nu_2(\th_{q,q+2i_2})-\(i-i_2\)
 \end{equation}
 hold.  \qed
 \end{theo}

 \begin{proof}
 We use again the notation $q_1=q/2$.
 To show the symmetry write
 $i=cQ+d$ with $Q-2q_1\le d\le Q$. In particular, we have $2q_1+d\ge
 Q$, and, hence, if $d\le Q-q_1$, the 2-adic
 valuation of
 $\th_{q,q+2i}$ is given by
 (cf.\ \eqref{ANa}, which also holds if $2q_1+d=Q$ as the
 quantity vanishes in this case)
 \begin{multline}
 \nu_2(\th_{q,q+2i})=\nu_2(\th_{q,q+2(cQ+d)})\\
 =-\sum _{j=0} ^{d-1}s(j)+\sum _{j=2q_1+d-Q} ^{2q_1+2d-Q-1}s(j)
 +(2q_1+d-Q)(s(c)-s(c+1)+1),
 \label{AU}
 \end{multline}
 while, if $d> Q-q_1$, it is given by
 (cf.\ \eqref{ANb})
 \begin{multline}
 \nu_2(\th_{q,q+2i})=\nu_2(\th_{q,q+2(cQ+d)})\\
 -\sum _{j=0} ^{d-1}s(j)+\sum _{j=2q_1+d-Q} ^{2q_1+2d-Q-1}s(j)
 -(2q_1+2d-2Q)+(Q-d)(s(c)-s(c+1)+1).
 \label{AV}
 \end{multline}
 Now let $i>i_2$, that is, $d>Q-q_1$. Symmetry around the center
 $i_2=(c+1)Q-q_1$ means to show that
 $$\nu_2(\th_{q,q+2(cQ+2Q-q-d)})=\nu_2(\th_{q,q+2(cQ+d)}).$$
 Using \eqref{AU}, the left-hand term is given by
 \begin{align*}
 \nu_2(&\th_{q,q+2((c+2)Q-2q_1-d)})\\
 &=
 -\sum _{j=0} ^{2Q-2q_1-d-1}s(j)
 +\sum _{j=Q-d} ^{3Q-2q_1-2d-1}s(j)
 +(Q-d)(s(c)-s(c+1)+1)\\
 &=-\sum _{j=0} ^{Q-1}s(j)+\sum _{j=2Q-2q_1-d} ^{Q-1}s(j)
 +\sum _{j=0} ^{2Q-1}s(j) -\sum _{j=0} ^{Q-d-1}s(j)
 -\sum _{j=3Q-2q_1-2d} ^{2Q-1}s(j)\\
 &\kern6cm
 +(Q-d)(s(c)-s(c+1)+1)\\
 &=\sum _{j=2Q-2q_1-d} ^{Q-1}s(j)
 +\sum _{j=Q} ^{2Q-1}s(j)
 -\sum _{j=0} ^{Q-d-1}s(j)
 -\sum _{j=3Q-2q_1-2d} ^{2Q-1}s(j)\\
 &\kern4cm
 +(Q-d)(s(c)-s(c+1)+1)\\
 &=\sum _{j=2Q-2q_1-d} ^{Q-1}s(j)
 +\sum _{j=Q-d} ^{Q-1}s(j)+Q
 -\sum _{j=3Q-2q_1-2d} ^{2Q-1}s(j)
 +(Q-d)(s(c)-s(c+1)+1).
 \end{align*}
 Now we apply the reflection identity \eqref{AI} to all the sum of digit
 functions. Thus, we obtain
 \begin{align*}
 \nu_2(&\th_{q,q+2((c+2)Q-2q_1-d)})
 =(2q_1+2d-Q)\log_2Q-\sum _{j=0} ^{2q_1+d-Q-1}s(j)
 -\sum _{j=0} ^{d-1}s(j)\\
 &\kern1cm
 -(2q_1+2d-Q)\log_2(2Q)+\sum _{j=0} ^{2q_1+2d-Q-1}s(j)
 +Q+(Q-d)(s(c)-s(c+1)+1)\\
 &=
 -\sum _{j=0} ^{d-1}s(j)
 +\sum _{j=2q_1+d-Q} ^{2q_1+2d-Q-1}s(j)
 -(2q_1+2d-2Q)+(Q-d)(s(c)-s(c+1)+1)\\
 &=\nu_2(\th_{q,q+2(cQ+d)}),
 \end{align*}
 by \eqref{AV}, as desired.
 That the values at the extreme points of the interval $i=(c+1)Q-q$ and
 $i=(c+1)Q$ are $0$, was already shown in Theorem~\ref{thm:2}, \textsc{Case~1}.

 Next we determine the 2-adic value of $\th_{q,q+2i}$ at the
 center
 $i_2=(c+1)Q-q/2$ of the interval. By \eqref{AU}, we have
 \begin{align*} \nu_2(\th_{q,q+2i_2})&=
 \nu_2(\th_{q,q+2((c+1)Q-q_1)})\\
 &=
 -\sum _{j=0} ^{Q-q_1-1}s(j)
 +\sum _{j=q_1} ^{Q-1}s(j)
 +q_1(s(c)-s(c+1)+1)\\
 &=
 -\sum _{j=0} ^{q_1-1}s(j)
 +\sum _{j=Q-q_1} ^{Q-1}s(j)
 +q_1(s(c)-s(c+1)+1)\\
 &=
 q_1\log_2Q
 -2\sum _{j=0} ^{q_1-1}s(j)
 +q_1(s(c)-s(c+1)+1).
 \end{align*}

 Inequality \eqref{AT} was already implicitly proved in the proof of
 Theorem~\ref{thm:2}, \textsc{Case~1d}. Namely, if $i=cQ+d$,
 with $d\ge Q-q_1$, then the conditions of \textsc{Case~1d} are satisfied, and
 there it was shown (see the paragraph after \eqref{ANb}) that
 $$\nu_2(\th_{q,q+2i})=\nu_2(\th_{q,q+2(cQ+d)})\ge Q-d=(c+1)Q-i,$$
 as desired.

 To prove inequality \eqref{AS}, we compute the difference of the 2-adic
 valuations of $\th_{q,q+2i_2}$ and $\th_{q,q+2i}$ for
 $i\in[i_2,(c+1)Q]$: let again $i=cQ+d$,
 $Q-q_1\le d\le Q$. Then, using \eqref{AV} again, we have
 \begin{align*}
 \nu_2&(\th_{q,q+2i_2})-\nu_2(\th_{q,q+2i})=
 \nu_2(\th_{q,q+2((c+1)Q-q_1)})-\nu_2(\th_{q,q+2(cQ+d)})\\
 &=-\sum _{j=0} ^{Q-q_1-1}s(j)
 +\sum _{j=q_1} ^{Q-1}s(j)
 +q_1(s(c)-s(c+1)+1)
 \\
 &\kern1cm
 +\sum _{j=0} ^{d-1}s(j)
 -\sum _{j=2q_1+d-Q} ^{2q_1+2d-Q-1}s(j)
 +(2q_1+2d-2Q)-(Q-d)(s(c)-s(c+1)+1)\\
 &=\sum _{j=Q-q_1} ^{d-1}s(j)
 +\sum _{j=q_1} ^{2q_1+d-Q-1}s(j)
 -\sum _{j=Q} ^{2q_1+2d-Q-1}s(j)
 +(q_1+d-Q)(s(c)-s(c+1)+3)\\
 &=\sum _{j=Q-q_1} ^{d-1}s(j)
 +\sum _{j=q_1} ^{2q_1+d-Q-1}s(j)
 -\sum _{j=0} ^{2q_1+2d-2Q-1}s(j)
 +(q_1+d-Q)(s(c)-s(c+1)+1)\\
 &=\sum _{j=Q-q_1} ^{d-1}s(j)
 +\sum _{j=q_1} ^{2q_1+d-Q-1}s(j)
 -2\sum _{j=0} ^{q_1+d-Q-1}s(j)
 +(q_1+d-Q)(s(c)-s(c+1)).
 \end{align*}
 Since $Q-q_1\ge q_1+d-Q$ and $q_1\ge q_1+d-Q$,
 we may apply
 the inequality (\ref{basins1}) twice to obtain
 $$
 \nu_2(\th_{q,q+2i_2})-\nu_2(\th_{q,q+2i})
 \ge(q_1+d-Q)(s(c)-s(c+1)+2).
 $$
 As we used already quite often, $s(c)-s(c+1)+1\ge0$. Therefore,
 $$\nu_2(\th_{q,q+2i_2})-\nu_2(\th_{q,q+2i})
 \ge q_1+d-Q=i-i_2,$$
 as desired.
\end{proof}

 \begin{corol} \label{cor:9}
 For a fixed even $q$, the values of the sequence
 $(\nu_2(\th_{q,q+2i}))_{i\ge0}$ are never $1$ for
 $i\in [(c+1)Q-q,(c+1)Q]$, except if $q$ is a power of $2$ and $c$ is
 even. In the
 latter case, $\nu_2(\th_{q,q+2i})$ is equal to $1$ for $i=(c+1)Q-q-1$
 and for $i=(c+1)Q-1$.

 \end{corol}

 \begin{proof} By the symmetry of the values around the center of the
 interval, and by the inequality \eqref{AT}, the only possible
 candidates are $i=(c+1)Q-q+1$ and $i=(c+1)Q-1$. Let us concentrate
 on $i=(c+1)Q-1$, which, as before, we write using the parameter $d$
 as  $i=(c+1)Q-1=cQ+d$.
 The formula \eqref{ANb} in
 the proof of Theorem~\ref{thm:2}, \textsc{Case~1d}, (which we used to prove
 \eqref{AT}), yields that $\nu_2(\th_{q,q+2i})$ is equal to
 \begin{align*}
 -\sum _{j=0} ^{Q-2}s(j)+\sum _{j=q-1} ^{q+Q-3}s(j)
 &-(q-2)+(s(c)-s(c+1)+1)\\
 &=-\sum _{j=0} ^{Q-1}s(j)+s(Q-1)+\sum _{j=q-1} ^{q+Q-2}s(j)-s(q+Q-2)\\
 &\kern5cm
 -(q-2)+(s(c)-s(c+1)+1)\\
 &=\cl{\log_2q}-s(q+Q-2)+1+(s(c)-s(c+1)+1).
 \end{align*}
 For this expression to be equal to $1$, we must have $s(q+Q-2)=\cl{\log_2q}$
 and $s(c)-s(c+1)+1=0$. The former is the case if and only if $q=Q$, that
 is, if $q$ is a power of $2$, and the latter is the case if and only
 if $c$ is even.
\end{proof}

 In an analogous manner, one can prove the following two theorems,
 with accompanying corollaries,
 covering the case where $q$ is odd.

 \begin{theo} \label{thm:10}
 Let $q$ be a fixed positive odd integer,
 $Q=2^{\cl{\log_2 q}}$, and let $c$ be any nonnegative integer.
 Then the values of the sequence
 $(\nu_2(\th_{q,q+2i}))_{i\ge 0}$ for
 $i\in [cQ,(c+1)Q-q]$ are symmetric around its center $cQ+(Q-q)/2$.
 The values at the extreme points of the interval
 $i=cQ$ and $i=(c+1)Q-q$ are $0$, and
 the value at the central points $i_3=cQ+(Q-q-1)/2$ and
 $cQ+(Q-q+1)/2$ is
 $$\nu_2(\th_{q,q+2i_3})=
 \frac {Q-q-1} {2}\log_2\frac {Q} {2}
 -2\sum _{j=0} ^{\frac {Q-q-3} {2}}s(j)
 -s\(\frac {Q-q-1} {2}\),
 $$
 independent of $c$.  Furthermore,
 for $i\in [cQ,i_3]$ the inequalities
 \begin{equation} \label{AW}
 \nu_2(\th_{q,q+2i})\ge i-cQ
 \end{equation}
 and
 \begin{equation} \label{AX}
 \nu_2(\th_{q,q+2i})\le \nu_2(\th_{q,q+2i_3})-\(i_3-i\)
 \end{equation}
 hold.
 \end{theo}

 \begin{corol} \label{cor:10}
 For a fixed odd $q$, the values of the sequence
 $(\nu_2(\th_{q,q+2i}))_{i\ge0}$ are never $1$ for
 $i\in [cQ,(c+1)Q-q]$, except if $q$ has the form $2^M+2^m-1$, for
 some positive integers $m$ and $M$, $m<M$. In the latter case, $\nu_2(\th_{q,q+2i})$
 is equal to $1$ for $i=cQ+1$ and for $i=(c+1)Q-q-1$.

 \end{corol}

 \begin{proof} The arguments from the proof of Corollary~\ref{cor:8}
 apply also here. Thus, again, the only possible
 candidates are $i=cQ+1$ and $i=(c+1)Q-q-1$. Furthermore, we must have
 $$s\(q+1-\frac {Q} {2}\)=1.$$
 This is only possible if $q+1-\frac {Q} {2}$ is a power of $2$,
 which means that $q$ has the form given in the statement of the corollary.
 \end{proof}

 \begin{theo} \label{thm:11}
 Let $q$ be a fixed positive odd integer,
 $Q=2^{\cl{\log_2 q}}$, and let $c$ be any nonnegative integer.
 Then the values of the sequence
 $(\nu_2(\th_{q,q+2i}))_{i\ge 0}$ for
 $i\in [(c+1)Q-q,(c+1)Q]$ are symmetric around its center $(c+1)Q-q/2$.
 The values at the extreme points of the interval $i=(c+1)Q-q$ and
 $i=(c+1)Q$ are $0$, and
 the value at the central points $i_4=(c+1)Q-(q+1)/2$ and
 $(c+1)Q-(q-1)/2$ is
 $$\nu_2(\th_{q,q+2i_4})=
 \frac {q-1} {2}\log_2Q
 -2\sum _{j=0} ^{\frac {q-3} {2}}s(j)
 -s\(\frac {q-1} {2}\)+\frac {q-1} {2}(s(c)-s(c+1)+1).
 $$
 Furthermore,
 for $i\in [i_4,(c+1)Q]$ the inequalities
 \begin{equation} \label{AY}
 \nu_2(\th_{q,q+2i})\ge (c+1)Q-i
 \end{equation}
 and
 \begin{equation} \label{AZ}
 \nu_2(\th_{q,q+2i})\le \nu_2(\th_{q,q+2i_4})-\(i-i_4\)
 \end{equation}
 hold.
 \end{theo}

 \begin{corol} \label{cor:11}
 For a fixed odd $q$, the values of the sequence
 $(\nu_2(\th_{q,q+2i}))_{i\ge0}$ are never $1$ for
 $i\in [(c+1)Q-q,(c+1)Q]$.

 \end{corol}

 \begin{proof}  The arguments from the proof of Corollary~\ref{cor:9}
 apply also here. The conclusion was that we can have
 $\nu_2(\th_{q,q+2i})=1$, for some $i$, only if $q$ is a power of
 $2$. This is a contradiction to our assumption that $q$ is odd.
 \end{proof}

 \section{Skew symmetric matrices
 and the parity of $\varepsilon_{2p,n}$}

 Let $\F=\C,\R$ and denote by $\gl(n,\F)\subset \rM_n(\F)$
 the group
 of $n\times n$ invertible matrices.  Recall
 that $\rA_n(\F)$ is the linear space of $n\times n$ skew symmetric matrices
 $A$ of order $n$
 over $\F$, i.e., $A\trans =-A$.  Clearly $\dim
 \rA_n(\F)={\binom n 2}$.
 Two matrices $A,B\in\rA_n(\F)$ are called
 \emph{congruent} if $A=TBT\trans$ for some $T\in \gl(n,\F)$.
 Let
 $S_2:=\left(\begin{smallmatrix}\hphantom{-}0&1\\-1&0\end{smallmatrix}\right)$.
 The following result is well-known in the real case, but its
 complex version
 does not seem to appear in standard modern books on linear algebra.

 \begin{prop}\label{eqvcsmat}  Let $\F=\R,\C$ and $A
 \in\rA_n(\F)$.  Then $A$ has
 even rank,
 $2p$ say, and $A$ is
 congruent over $\F$ to a direct sum of $p$ copies of $S_2$ and the
 $(n-2p)\times (n-2p)$ zero matrix.  In particular,
 $B\in \rA_n(\F)$ is congruent to $A$ over $\F$ if and only if
 $\rank A=\rank B$.
 \end{prop}
 \begin{proof}  We first prove the fact that any $A\in\rA_n(\F)$ is
 congruent to the direct sum of copies of $S_2$ and $0$.
 The result is trivial if $A=0$.
 Let $n=2$ and $\rank A=2$.  Then $A=aS_2$ for some $0\ne a\in\F$.
 For $\F=\C$ we have $A=(\sqrt{a} I_2)S_2 (\sqrt{a} I_2)\trans$.
 For $\F=\R$ and $a>0$ the above formula holds.  For $a<0$ we have
 $A=(\sqrt{-a} P)S_2 (\sqrt{-a} P)\trans$, where
 $P:=\left(\begin{smallmatrix}0&1\\1&0\end{smallmatrix}\right)$.

 Assume by induction that any $A\in\rA_n(\F)$ is congruent
 to the direct sum of
 copies of $S_2$ and $0\in\rM_n(\F)$ for $n=m\ge 2$.  Let $n=m+1$
 and $A\in \rA_{m+1}(\F)$.  Suppose first that $\det A=0$.  Let
 $\0\ne\x\in\F^n$ and $A\x=0$.  Let $Q\in\gl(m+1,F)$ such that the
 last column of $Q\trans$ is $\x$.  Then $QAQ\trans =A_1\oplus 0$,
 and $A_1\in \rA_{m}(\F)$.  Use the induction hypothesis to deduce
 that $A$ is conjugate to the direct sum of copies of $S_2$ and
 $0$.
 It remains to study the case where $m+1$ is even and $\det
 A\ne 0$.  Let $A=(a_{ij})_{i,j=1}^{m+1}$ and
 $A_1=(a_{ij})_{i,j=1}^m\in \rA_m(\F)$.  Since $A$ has rank $m+1$,
 $A_1$ has at least rank $m-1$.  Since $m$ is odd $A_1$ has
 exactly rank $m-1$.  So $A_1$ is conjugate to a direct sum of
 $\frac{m-1}{2}$ copies of $S_2$ and one copy of $0$.
 Using the corresponding
 congruence on $A$, we may assume without loss of
 generality that $A_1=
 \left(\bigoplus_{i=1}^{(m-1)/2} S_2\right) \oplus 0 $.
 For each
 $i=1,\ldots,m-1$ subtract from
 column $m+1$ of
 $A$ the corresponding multiple of column $i$ to obtain
 the zero
 element
 for the $(i,m+1)$-entry.  Repeat these elementary operations
 with the rows of $A$ to eliminate the
 $(m+1,i)$-entry for
 $i=1,\ldots,m-1$.  The
 resulting matrix is of the form
 $B= \left(\bigoplus_{i=1}^{(m-1/2)} S_2\right)\oplus bS_2$.
 Clearly $A$ is congruent to $B$.  Hence $b\ne
 0$. Since $bS_2$ is congruent
 to $S_2$, we deduce that $A$ is congruent to
 a direct sum of
 copies of $S_2$.

 Since a direct sum of copies of $S_2$ and $0$ has an even rank
 we deduce that any $A\in\rA_n(\F)$ has
 even rank.
\end{proof}

 The following result is known to the experts. We bring its
 proof for completeness.
 \begin{prop}\label{irw2pc}  Let $\F=\C,\R$, $n\ge 2$, $p\ge 1$ be integers and
 assume that $p\le \lfloor \frac{n}{2}\rfloor$.
 Let $\P \rW_{2p,n}(\F)\subset \P \rA_n(\F)$ be the projective variety of all
 {\em(}nonzero{\em)} skew symmetric matrices of rank
 at most $2p$. Then $\P
 \rW_{2p,n}(\F)$
 is an irreducible projective variety in
 $\P \rA_n(\F)$ of codimension $\binom {n-2p} 2$.
 The variety of its singular points is $\P \rW_{2(p-1),n}(\F)$.
 \end{prop}

 \begin{proof}  Let $\rW_{2p,n}^o(\F)\subset \rA_n(\F)$ be the quasi-variety of all
 $n\times n$ skew symmetric matrices of rank $2p$.
 Proposition~\ref{eqvcsmat} yields that $\gl(n,\F)$ acts transitively on
 $\rW_{2p,n}^o(\F)$.  Hence $\rW_{2p,n}^o(\F)$ is a homogeneous
 space and a manifold.  Since $\rW_{2(p-1),n}(\F)$ is a strict affine
 subvariety of $\rW_{2p,n}(\F)$ it follows that $\rW_{2p,n}^o(\F)$
 is a subset of smooth points and $\dim \rW_{2p,n}^o(\F)=\dim
 \rW_{2p,n}(\F)$.  The neighborhood of each point is obtained by the corresponding
 action of the neighborhood of $I_n\in \gl(n,\F)$.
% In particular, the tangent bundle of $A\in\rW_{2p,n}^o(\F)$ is the quotient of
% $\rA_n(\F)$ by the subspace $\{X\in \rA_n(\F):\;XA+AX\trans=0\}$.
 Proposition~\ref{eqvcsmat} yields that the orbit of any $B\in \rW_{2(p-1),n}(\F)$
 does not contain any matrix in $\rW_{2p,n}^o(\F)$.  Hence $\rW_{2(p-1),n}(\F)$ is the
 variety of singular points in $\rW_{2p,n}(\F)$.

 We now find the dimension and codimension of $\rW_{2p,n}(\F)$.  For $p=\lfloor
 \frac{n}{2}\rfloor$, we have $\rW_{2p,n}(\F)=\rA_n(\F)$.  Hence $\dim \rW_{2\lfloor
 \frac{n}{2}\rfloor,n}(\F)= {\binom n 2}$ and $\codim \rW_{2\lfloor
 \frac{n}{2}\rfloor,n}(\F)= \binom {n -2\lfloor\frac{n}{2}\rfloor }
 2=0$.  Let $1\le p < \lfloor\frac{n}{2}\rfloor$.  Let $A=(a_{ij})_{i,j=1}^n\in
 \rW_{2p,n}^o(\F)$.  Then $A$ has $2p$ independent rows.  Assume
 for simplicity
 that the first $2p$ rows are
 linearly independent.  Hence the first $2p$ columns of $A$ are
 linearly independent.  Hence $A_1=(a_{ij})_{i,j=1}^{2p}\in
 \rA_{2p}(\F)$ is nonsingular.  Therefore there exists a unique block
 lower triangular matrix
 $T=\left(\begin{smallmatrix} I_{2p} &0\\R& I_{n-2p}\end{smallmatrix}\right)$
 such that $TAT\trans =A_1\oplus 0$.
 Equivalently $A=T^{-1}(A_1\oplus 0)(T\trans)^{-1}$.
 Hence $\dim \rW_{2p,n}=\binom {2p} 2 + 2p(n-2p)$ and
 $\codim \rW_{2p,n}=\binom {n-2p}2$.
 Thus, $\codim \P\rW_{2p,n}=\binom {n-2p} 2$.
\end{proof}

 \begin{theo}\label{degskn}  Let $\F=\R,\C$ and $4 \le n, 1\le p <
 \lfloor\frac{n}{2}\rfloor$.  Let $\P\rW_{2p,n}(\F)$ be the
 irreducible variety of all {\em(}projectivized{\em)} nonzero skew symmetric
 $n\times n$ matrices of rank
 at most $2p$ in the projective space
 $\P\rA_n(\F)$ of all nonzero $n\times n$ skew symmetric matrices
 over $\F$.  Then the degree of $\P\rW_{2p,n}(\C)$ is odd if and
 only if either $p$ or $n-p$ is divisible by $2^{\lceil \log_2 (n-2p)\rceil}$.
 Furthermore, if either $p$ or $n-p$ is divisible by $2^{\lceil \log_2
 (n-2p)\rceil}$, then any $\left(\binom {n-p} 2+1\right)$-dimensional subspace of $n\times n$
 real skew symmetric matrices
 contains a nonzero matrix of rank
 at most $2p$. For these values
 of $n$ and $p$ the dimensions of subspaces are best possible, i.e.,
 \eqref{asrc} holds.
 \end{theo}
 \begin{proof}  Recall from (\ref{degvk}) that $\varepsilon_{2p,n}:= \deg \P \rW_{2p,n}
 =\delta_{2p+1,n}/2^{n-2p-1}$.  The definition
 (\ref{deftqn}) of $\th_{q,n}$ yields that
 $\nu_2(\varepsilon_{2p,n})=\nu_2(\th_{n-2p-1,n})-(n-2p-1)$.
 Proposition~\ref{prop:1} yields that $\nu_2(\varepsilon_{2p,n})=
 \nu_2(\th_{n-2p,n})$.  Use Theorem~\ref{thm:2} to deduce that
 $\varepsilon_{2p,n}$ is odd if and only if
 either $p$ or $n-p$ is divisible by $2^{\lceil \log_2 (n-2p)\rceil}$.

 Assume that  either $p$ or $n-p$ is divisible by $2^{\lceil \log_2 (n-2p)\rceil}$.
 Then the discussion in \S1 implies that any $\left(\binom {n-p}
 2+1\right)$-dimensional subspace of $n\times n$ real skew symmetric matrices
 contains a nonzero matrix of rank
 at most $2p$.  The sharpness
 of these dimensions follows from the fact that a complex subspace
 $L$ of $\rA_n(\C)$ of dimension $\binom {n-p} 2$ in general
 position will not contain a nonzero $A\in \rA_n(\C)$ of rank
 at most $2p$.
\end{proof}

 \begin{corol}\label{subcsm}
 Let $n\equiv 2\;({\rm mod\;} 4)$.  Then any two-dimensional real
 subspace of $n\times n$ skew symmetric matrices contains a
 nonzero singular matrix.
 \end{corol}

 \section{Rectangular matrices
 and the parity of $\gamma_{k,m,n}$}

 In this section we consider the parity problem for
 $\gamma_{k,m,n}$, the latter being defined in (\ref{degkmn}).  It is more convenient
 to introduce the following symmetric quantity.
 For $n\in\N$, let $\h(n)$ be the \emph{hyperfactorial} $\prod_{k=0}^{n-1} k!$.
 Let $a,b,c\in\N$.  Then a straightforward calculation shows:
 \begin{align} \nonumber
 B(a,b,c)&=\prod_{i=1}^a\prod_{j=1}^b\prod_{k=1}^c
 \frac {(i+j+k-1)} {(i+j+k-2)}\\
 \nonumber
 &=\prod _{i=1} ^{a}\frac {(b+c+i-1)!\,(i-1)!} {(b+i-1)!\,(c+i-1)!}\\
  \label{eq:M2}
 &=\frac {\h(a)\,\h(b)\,\h(c)\,\h(a+b+c)}
 {\h(a+b)\,\h(b+c)\,\h(c+a)}.
 \end{align}
 $B(a,b,c)$ is a symmetric function on $\N^3$.
 We remark that $B(a,b,c)$ is
 the number of plane partitions which are contained in an
 $a\times b\times c$ box (see e.g.\
 \cite{Bressoud}).
 From the definition (\ref{degkmn}) of $\gamma_{k,m,n}$ it is obvious that
 \begin{equation}\label{gambeq}
 \gamma_{k,m,n}= B(n-k,m-k,k), \quad 1\le k\le \min(m,n).
 \end{equation}

%June 25, 2005

 For simplicity of notation we let
 $$S(a):=\sum_{i=0}^{a-1} s(i)\quad \text{and}\quad  \nu(a,b,c):=\nu_2(B(a,b,c))
 \quad \text{for any } a,b,c\in\N.$$
 Then,
 by Proposition~\ref{2adn!}, we have $\nu_2(\h(a))=\frac{(a-1)a}{2}-S(a)$ and
 \begin{equation}
 \label{eq:nu2}
 \nu(a,b,c)=-S(a+b+c)-S(a)-S(b)-S(c)+S(a+b)+S(b+c)+S(a+c).
 \end{equation}

 \begin{lemma}\label{basnuid}  Let $a,b,c\in\N$.  Then the
 following identities hold:
 \begin{align}
 \nu(2a,2b,2c)&=2\nu(a,b,c), \label{basnuid1}\\
 \nu(2a,2b+1,2c+1)&=\nu(a,b+1,c)+\nu(a,b,c+1),
 \label{basnuid2}\\
 \nu(2a+1,2b,2c)&=\nu(a,b,c)+\nu(a+1,b,c), \label{basnuid3}\\
 \nu(2a+1,2b+1,2c+1)&=\nu(a,b+1,c+1)+\nu(a+1,b,c) \nonumber\\
 &\kern2cm +s(b+c)-s(b+c+1)+2. \label{basnuid4}
 \end{align}
 In particular,
 \begin{enumerate}
 \item[\em(1)]
 $B(2a,2b,2c)$ is odd if and only if $B(a,b,c)$ is odd.
 \item[\em(2)]
 $B(2a,2b+1,2c+1)$ is odd if and only if both $B(a,b+1,c)$
 and $B(a,b,c+1)$ are odd.
 \item[\em(3)]
 $B(2a+1,2b,2c)$ is odd if and only if both $B(a,b,c)$
 and $B(a+1,b,c)$ are odd.
 \item[\em(4)]
 $B(2a+1,2b+1,2c+1)$ is always even.
 \end{enumerate}
 \end{lemma}

 \begin{proof} Equality (\ref{fundeq}) is equivalent to
 \begin{equation}\label{fundeqS}
 S(2p)=2S(p)+p,\;S(2p+1)=S(p+1)+S(p)+p, \quad \text{for any } p\in\N.
 \end{equation}
 Use (\ref{eq:nu2}) and the above equalities to deduce
 (\ref{basnuid1})--(\ref{basnuid3}) straightforwardly.
 In particular, these equalities yield the
 corresponding claims about the oddness of $B(u,v,w)$.

 To obtain  (\ref{basnuid4}), we use in addition to
 (\ref{eq:nu2}) and the above equalities the obvious equality
 $S(p+1)=S(p)+s(p)$ for $p=b+c,b+c+1$.  Since $s(p+1)\le s(p)+1$
 it follows that  $B(2a+1,2b+1,2c+1)$ is always even.
\end{proof}

 \begin{rem} \label{rem2}
 (1) Lemma~\ref{basnuid} provides us with an algorithm to compute
 the parity of $B(a,b,c)$ from the binary expansions of $a,b,c$
 directly, without having to actually compute $B(a,b,c)$.
 Namely, given $a,b,c$, one determines the parities of $a,b,c$.
 If all of $a,b,c$ are odd, then Conclusion~(4) in Lemma~\ref{basnuid}
 says that $B(a,b,c)$ is even. If all of $a,b,c$ are even, then one
 uses Conclusion~(1) to reduce the problem to the problem of
 determining the parity of $B(a/2,b/2,c/2)$. If exactly two of $a,b,c$ should be
 odd, then one uses Conclusion~(2) for a similar reduction, and
 if only one of $a,b,c$ is odd, then one uses Conclusion~(3).
 One sees quickly that this yields an algorithm which can be
 most conveniently run on the binary expansions of $a,b,c$.
 See Proposition~\ref{oddconb} for an attempt to turn this algorithm
 into a concrete characterization of those $a,b,c$ for which
 $B(a,b,c)$ is odd.

 (2) For the interested reader, we remark that the inspiration for
 Lemma~\ref{basnuid} comes from results on plane partitions due to
 Stembridge and Eisenk\"olbl. More precisely, Stembridge showed in
 \cite{Ste} that a certain $(-1)$-enumeration (for our purposes it
 suffices to say that this means a weighted enumeration in which
 some plane partitions count as $1$, as in ordinary enumeration,
 and others count as $-1$) of plane partitions contained in an
 $a\times b\times c$-box is equal (up to sign) to the number of
 self-complementary plane partitions contained in the same box.
 Since, by definition of self-complementary plane partitions, there
 cannot exist any if all of $a,b,c$ are odd, this result implies
 immediately Conclusion~(4) in Lemma~\ref{basnuid}. Subsequently, Eisenk\"olbl
 \cite{EisTAF} has embarked on the $(-1)$-enumeration of self-complementary
 plane partitions. Her result is that a certain $(-1)$-enumeration
 of self-complementary plane partitions contained in an
 $a\times b\times c$-box is equal to
 \begin{align*}
 B\(\tfrac a2,\tfrac b2, \tfrac c2\)^2 \quad &\text{for $a,b,c$ even,}\\
 B\(\tfrac a2,\tfrac{b+1}2,\tfrac{c-1}2\) B\(\tfrac
 a2,\tfrac{b-1}2,\tfrac{c+1}2\) \quad &\text{for $a$ even and $b$, $c$
   odd,}\\
 B\(\tfrac {a+1}2,\tfrac{b}2,\tfrac{c}2\) B\(\tfrac
 {a-1}2,\tfrac{b}2,\tfrac{c}2\)
 \quad &\text{for $a$ odd and $b$, $c$
   even.}
 \end{align*}
 The first of the three cases implies Conclusion~(1) in
 Lemma~\ref{basnuid}, the second implies Conclusion~(2),
 and the third implies Conclusion~(3). The relations
 (\ref{basnuid1})--(\ref{basnuid4}) refine these conclusions on the
 level of $2$-adic valuations.
 \end{rem}

 Let $b$ and $c$ be two nonnegative integers.
 Furthermore, let
 $$b=\sum_{i=0}^{\infty} b_i2^i, \; c=\sum_{i=0}^{\infty}c_i2^i, \quad
 b_i,c_i\in \{0,1\},\;i=0,\ldots,$$
 be the binary expansion of $b$ and $c$, respectively.
 We say that the pair $(b,c)$ has a \emph{disjoint binary expansion}
 if $b_ic_i=0$ for $i=0,\ldots$.
 The following proposition is straightforward to establish.
 \begin{prop} \label{no1com}  Let $b$ and $c$ be nonnegative integers.
 Then $s(b+c)\le s(b)+s(c)$, and $s(b+c)=s(b)+s(c)$ if and only
 the pair $(b,c)$ has a disjoint binary expansion.
 \end{prop}

 To find out
 under which conditions on the parameters $a,b,c$
 the number $B(a,b,c)$  is odd, it is enough to consider the
 case $a\le \min(b,c)$.
 \begin{theo}\label{a=123}  Let $b,c\in\N$.  Then
 \begin{enumerate}
 \item
 $B(1,b,c)$ is odd
 if and only if $(b,c)$ has a disjoint binary expansion.
 In particular, for any
 $q\ge 0$, $B(2^q,2^q b,2^q c)$ is odd
 if and only if $(2^q b,2^q c)$ has a disjoint binary expansion.

 \item Let $\min(b,c)\ge 2$.
 \begin{itemize}
 \item
 If $b$ and $c$ are even then $B(2,b,c)$ is odd if and only if
 the pair $(b,c)$ has a disjoint binary expansion.
 \item
 If $b$ is even and $c$ is odd
 then $B(2,b,c)$ is odd if and only if the pairs $(b,c)$ and $(b,c+1)$
 have disjoint binary expansions.
 \item
 If $b$ and $c$ are odd then
 $B(2,b,c)$ is odd if and only if the pairs $(b,c+1)$ and $(b+1,c)$
 have disjoint binary expansions.
 \end{itemize}
 \item
 Let $\min(b,c)\ge 3$.
 \begin{itemize}
 \item
 Assume that $b$ and $c$ are even.

 If $b,c\equiv 0$ \textrm{\em mod} $4$ then $B(3,b,c)$ is
 odd if and only if the pair $(b,c)$ has a disjoint binary expansion.

 If $b,c+2\equiv 0$ \textrm{\em mod} $4$ then $B(3,b,c)$ is
 odd if and only if the pairs $(b,c)$ and $(b,c+2)$ have disjoint binary expansions.

 If $b,c\equiv 2$ \textrm{\em mod} $4$ then
 $B(3,b,c)$ is
 odd if and only if the pairs $(b,c)$, $(b,c+2)$ and $(b+2,c)$ have disjoint
 binary expansions.
 \item
 Assume that $b$ is even and $c$ is odd.

 If $b\equiv 0$ \textrm{\em mod} $4$ then $B(3,b,c)$ is
 odd if and only if the pairs $(b,c)$ and $(b,c+1)$ have
 disjoint binary expansions.

 If $b\equiv 2$ \textrm{\em mod} $4$ then $B(3,b,c)$ is
 odd if and only if the pairs $(b,c+1)$ and $(b+2,c)$
 have disjoint binary expansions.
 \end{itemize}
 \end{enumerate}
 \end{theo}
 \begin{proof}
 \begin{enumerate}
 \item
 The expression (\ref{eq:M2}) yields that $B(1,b,c)=\binom {b+c} b$.
 Hence $\nu(1,b,c)=s(b)+s(c)-s(b+c)$.  Thus $\nu(1,b,c)=0$ if and
 only if $(b,c)$ has a disjoint binary expansion.
 The last
 assertion follows from (\ref{basnuid1}) and from the observation
 that $(b,c)$ has a disjoint binary expansion if and only if $(2^q b,2^q c)$
 has a disjoint binary expansion.

 \item
 \begin{itemize}
 \item
 If $b$ and $c$ are even then
 Conclusion~(1) in Lemma~\ref{basnuid} and item~1, which we just established, yield that
 $B(2,b,c)$ is odd if and only the pair $(b,c)$ has
 a disjoint binary expansion.
 \item
 Assume that $b=2b'$, $c=2c'+1$, $b',c'\in\N$.
 Since $\nu(a,b,c)$ is a symmetric function in $a,b,c$,
 (\ref{basnuid3}) yields that
 $\nu(2,2b',2c'+1)=\nu(1,b',c')+\nu(1,b',c'+1)$.
 Hence $B(2,b,c)$ is odd if and only if $(b',c')$ and $(b',c'+1)$
 have disjoint binary expansions.  This is equivalent to the assumption
 that the two even pairs $(b,c-1)$ and $(b,c+1)$ have disjoint binary
 expansions.  Since $b$ is even and $c$ odd the assumption that
 $(b,c-1)$ has disjoint binary
 expansions is equivalent to the assumption
 that $(b,c)$ has a disjoint binary
 expansion.
 \item
 Assume that $b=2b'+1$, $c=2c'+1$, $b',c'\in\N$.  Then
 (\ref{basnuid2}) yields
 $\nu(2,b,c)=\nu(1,b'+1,c)+\nu(1,b',c'+1)$.
 Hence $B(2,b,c)$ is odd if and only if $(b'+1,c')$ and $(b',c'+1)$
 have disjoint binary expansions.  This is equivalent
 to the assumption
 that the two even pairs $(b+1,c-1)$ and $(b-1,c+1)$ have disjoint binary
 expansions.  This is also equivalent to to the assumption
 that $(b+1,c)$ and $(b,c+1)$ have disjoint binary
 expansions.
 \end{itemize}

 \item
 \begin{itemize}
 \item
 Assume that $b=2b'$, $c=2c'$, $a',b'\ge 2$.
 Relation~(\ref{basnuid3}) implies that $B(3,b,c)=B(1,b',c')+B(2,b',c')$.
 Hence $B(3,b,c)$ is odd if an only if $B(1,b',c')$ and $B(2,b',c')$
 are odd.  Recall that $B(1,b',c')$ is odd if and only if
 $(b',c')$ has a disjoint binary expansion.  This is equivalent
 to the assumption that $(b,c)$ has a disjoint binary expansion.

 Assume that $b'$ and $c'$ are even.  Then $B(2,a',b')$ is
 odd if and only $(b',c')$ has a disjoint binary expansion.

 Assume now that $b'$ is even and $c'$ are is odd.
 Then $B(2,b',c')$ is odd if $(b',c')$ and $(b',c'+1)$
 have disjoint binary expansions.  This is equivalent to
 the assumption that $(b,c)$ and $(b,c+2)$ have
 disjoint binary expansions.

 Assume now that $b'$ and $c'$ is odd.  Then $B(2,b',c')$ is odd
 if and only if the pairs $(b',c'+1)$ and $(b'+1,c')$ have
 disjoint binary expansions.  This is equivalent to the
 assumption that the pairs $(b,c+2)$ and $(b+2,c)$
 have disjoint binary expansions.

 \item
 Assume that $b=2b'$, $c=2c'+1$, $b'\ge 2$, $c'\ge 1$.
 Then $\nu(3,2b',2c'+1)=\nu(1,b',c'+1)+\nu(2,b',c')$.
 Hence $B(3,b,c)$ is odd if and only if
 $B(1,b',c'+1)$ and $B(2,b',c')$ are odd.
 $B(1,b',c'+1)$ is odd if and only if $(b',c'+1)$
 has a disjoint binary expansion.  This is equivalent
 to the assumption that $(b,c+1)$ has a disjoint binary expansion.

 Assume that $b'$ is even, i.e., $4\mid b$.
 Suppose first that $c'$ is even, i.e., $4\mid (c-1)$.
 Then $B(2,b',c')$ is odd if and only if $(b',c')$
 has a disjoint binary expansion.  This is equivalent
 to the assumption that $(b,c)$ has a disjoint binary expansion.

 Assume second that $c'$ is odd, i.e., $4\mid (c+1)$.
 Then $B(2,b',c')$ is odd if and only if $(b',c')$ and
 $(b',c'+1)$ have disjoint binary expansions.
 This is equivalent to the assumption that
 $(b,c)$ and $(b,c+1)$ have disjoint binary expansions.

 Assume now that $b'$ is odd, i.e., $4\mid (b+2)$.
 Suppose first that $c'$ is even.  Then $B(2,b',c')$
 is odd if and only if $(b',c')$ and $(b'+1,c')$
 have disjoint binary expansions.
 This is equivalent to the assumption that $(b,c)$ and $(b+2, c)$
 have disjoint binary expansions.

 Suppose second that $c'$ is odd.
 Assume first that $c'=1$, i.e., $c=3$. Then $B(2,b',1)$ is odd if and only if
 $(b',2)$ has a disjoint binary expansion.  This is equivalent
 to the assumption that $(b,4)=(b,c+1)$ has a disjoint binary expansion.
 Since $4\mid (b+2)$ it follows that $(b+2,3)$ has a disjoint binary expansion.

 Assume second that $c'\ge 3$. Then $B(2,b',c')$ is odd if and only if $(b'+1,c')$
 and $(b',c'+1)$ have disjoint binary expansions.  This is
 equivalent to the assumption that $(b+2,c)$ and $(b,c+1)$ have
 disjoint binary expansions.
 \end{itemize}
 \end{enumerate}
 \end{proof}

 The above theorem can be generalized schematically as follows:

 \begin{prop}\label{oddconb}
 Let $a,b,c\in \N$ and assume that $\min(b,c)\ge a$.  Let
 $q:=\lceil \log_2 a\rceil$ and assume that $b\equiv b_r$, $c\equiv
 c_r$ \textrm{\em mod} $2^q$ for some
 $b_r,c_r\in [0,2^q-1]$. Then there exists a sequence of nonnegative
 integers $d_i,e_i \in [0,2^q-1]$, $i=1,\ldots, N(a,b_r,c_r)$, depending only on
 $a,b_r,c_r$, such
 that $B(a,b,c)$ is odd if and only $(b+d_i,c+e_i)$ has a
 disjoint binary expansion for $i=1,\ldots,N(a,b_r,c_r)$.
 \end{prop}

 \begin{proof}
 We prove the proposition by induction on $q$.
 For $q=0,1$ the proposition holds in view of Theorem~\ref{a=123}.
 Assume that the proposition holds for any $q\le p-1$, where $p\ge
 2$ and any $b,c$ such that $\min(b,c)\ge a$.  Assume that
 $\lceil \log_2 a\rceil=p$.

 \begin{itemize}
 \item
 Let $a=2a'$, $a'\in\N$.
 Then $\lceil \log_2 a\rceil= \lceil \log_2 a'\rceil + 1$.
 Assume first that $b$ and $c$ are even.
 Then $\nu(a,b,c)=2\nu(\frac{a}{2},\frac{b}{2},\frac{c}{2})$ and
 the proposition follows straightforwardly from the induction hypothesis.

 Assume now that $b=2b'$, $c=2c'+1$, $a'\le b',c'\in\N$.
 Then $\nu(2a',2b',2c'+1)=\nu(a',b',c')+\nu(a',b',c'+1)$.
 Hence $B(a,b,c)$ is odd if and only $B(a',b',c')$ and
 $B(a',b',c'+1)$ are odd.  Use the induction hypothesis and the
 observation that $(u,v)$ has a disjoint binary expansion
 if and only if $(2u,2v+1)$ has one for the case $B(a',b',c')$
 to deduce the proposition.  (For the case $B(a',b',c'+1)$
 note that if $d_i\le 2^{p-1}-1$
 then $2d_i+1 \le 2^p-1$.)

 Use the equality $B(a,b,c)=B(a,c,b)$ to deduce the proposition
 in the case that $b$ is odd and $c$ is even.

 Assume now that $b=2b'+1$, $c=2c'+1$, $a'\le b',c'\in\N$.
 Then $\nu(2a',2b'+1,2c'+1)=\nu(a',b'+1,c')+\nu(a',b',c'+1)$.
 Hence $B(a,b,c)$ is odd if and only $B(a',b'+1,c')$ and
 $B(a',b',c'+1)$ are odd.  Use the induction hypothesis
 and the above remarks to deduce the proposition in this case.

 \item
 Assume that $a=2a'+1$, $a'\in \N$.  Then $\lceil \log_2 a\rceil
 =\lceil \log_2 a+1\rceil=
 \lceil \log_2 a'+1\rceil + 1$.

 Since $B(a,b,c)=B(a,c,b)$ is even if $a,b,c$ are odd,
 it is enough to consider the case $b=2b'$, $a'<b'\in\N$.
 Assume first that $c=2c'$, $a'<c'\in\N$.
 Then $\nu(2a'+1,2b',2c')=\nu(a',b',c')+\nu(a'+1,b',c')$.
 Hence $B(a,b,c)$ is odd if and only $B(a',b',c')$ and
 $B(a'+1,b',c')$ are odd.
 Use the induction hypothesis and the above arguments to deduce
 the proposition.

 Assume finally that $c=2c'+1$, $a'\le c'$.
 Then $\nu(2a'+1,2b',2c'+1)=\nu(a',b',c'+1)+\nu(a'+1,b',c')$.
 Hence $B(a,b,c)$ is odd if and only $B(a',b',c'+1)$ and
 $B(a'+1,b',c')$ are odd.  The case where $B(a',b',c'+1)$ is odd
 is done by induction on $a'$.  For $c'>a'$ the case
 $B(a'+1,b'c')$ is odd is done by induction on $a'+1$.
 The case $c'=a'$, which is equivalent to the case
 where $B(a'+1,b',a')=B(a',b',a'+1)$ is odd is done by induction
 on $a'$.  (As in the case of $B(3,2b',3)$ in Theorem~\ref{a=123}.)
 \end{itemize}
 \end{proof}

 In principle the proof of Proposition~\ref{oddconb} can be used
 to find the sequences $d_i,c_i$, $i=1,\ldots,N(a,b_r,c_r)$,
 recursively.  However, the explicit construction of all such sequences seems
 complicated even in the simple case where $a=2^q$, $q=2,\ldots$.
 Note that Case~1 of Theorem~\ref{a=123} finds the sequence
 for $a=2^q$ and $b_r=c_r=0$.  The cases $b_r=0$, $c_r=1$ and
 $b_r=c_r=1$ have simple results.

 \begin{prop}\label{a=2qbc=1}
 Let $q\in\N$ and $2^q\le b,c\in\N$.  If $2^q\mid b$, $2^q\mid (c-1)$, then
 $B(2^q,b,c)$ is odd if and only if $(b,c)$ and $(b,c+2^{q}-1)$
 have disjoint binary expansions.
 If $2^q\mid (b-1)$, $2^q\mid (c-1)$, then
 $B(2^q,b,c)$ is odd if and only if $(b+1,c)$, $(b,c+1)$, $(b+2^{q}-1,c)$, and $(b,c+2^{q}-1)$
 have disjoint binary expansions.
 \end{prop}

 The proof of this proposition is left to the reader.

 The results of \S1 yield the following theorem.
 \begin{theo}\label{dkmncr}  Let $k,m,n\in\N$, assume that
 $k< \min(m,n)$, and let $\gamma_{k,m,n}$ be the positive integer given
 by \eqref{degkmn}.
 Let $\rM_{m,n}(\R)$ be the space of all
 $m\times n$ real valued matrices.  Let $L\subset \rM_{m,n}(\R)$
 be a subspace of dimension $(m-k)(n-k)+1$.  If $\gamma_{k,m,n}$
 is odd then $L$ contains a nonzero matrix rank
 at most $k$.
 \end{theo}

 \begin{corol}\label{cordkmn} For the following positive integers
 $1\le k <n\le m$ any $((m-k)(n-k)+1)$-dimensional subspace of
 $\rM_{m,n}(\R)$ contains a nonzero matrix rank
 at most $k$:
 \begin{enumerate}
 \item
 $k=n-1$ and $(m-n+1,n-1)$ has a disjoint binary expansion.
 \item
  $2 \le k=n-2$.
  \begin{itemize}
  \item
  $n$ and $m$ are even and $(n-2,m-n+2)$ has a disjoint binary expansion.
  \item
  $n$ is even, $m$ is odd, $(n-2,m-n+2)$ and $(n-2,m-n+3)$
  have disjoint binary expansions.
  \item
  $n$ is odd, $m$ is even, $(n-2,m-n+3)$ and $( n-1,m-n+2)$
  have disjoint binary expansions.
  \end{itemize}

  \item
  $3\le k=n-3$.
  \begin{itemize}
  \item
  $4\mid (n-3)$, $4\mid m$, and $(n-3,m-n+3)$ has a disjoint binary
  expansion.
  \item
  $4\mid (n-3)$, $4\mid (m+2)$, $(n-3,m-n+3)$ and $(n-3,m-n+5)$
  have disjoint binary expansions.
  \item
  $4\mid (n-1)$, $4\mid (m+2)$, $(n-3,m-n+3)$ and $(n-1,m-n+3)$
  have disjoint binary expansions.
  \item
  $4\mid (n-1)$, $4\mid m$, $(n-1,m-n+3)$ and $(n-3, m-n+5)$
  have disjoint binary expansions.
  \item
  $4\mid (n-3)$, $m$ odd, $(n-3,m-n+3)$ and $(n-3, m-n+4)$
  have disjoint binary expansions.
  \item
  $4\mid (n-1)$, $m$ odd, $(n-3,m-n+4)$ and $(n-1,m-n+3)$
  have disjoint binary expansions.
  \item
  $4\mid (m-n+3)$, $n$ is even, $(n-3,m-n+3)$ and $(n-2, m-n+3)$
  have disjoint binary expansions.
  \item
  $4\mid (m-n+5)$, $n$ is even, $(n-2,m-n+3)$ and $(n-3, m-n+5)$
  have disjoint binary expansions.

  \end{itemize}

 \item
 Let $q\in\N$.
 \begin{itemize}
 \item
  $n=k+2^q$, $2^q\mid k$, $2^q\mid m$, and $(k,m-k)$
  has a disjoint binary expansion.
  \item
  $n=k+2^q$, $2^q\mid k$, $2^q\mid (m-1)$, $(k,m-k)$ and $(k,m-k+2^q-1)$
  have disjoint binary expansions.
  \item
  $2^{q+1}<n=k+2^q$, $2^q\mid (k-1)$, $2^q\mid (m-1)$, $(k,m-k)$ and
  $(k+2^q-1,m-k)$
  have disjoint binary expansions.
  \item
  $2^{q+1}<n=k+2^q$, $2^q\mid (k-1)$, $2^q\mid (m-2)$, and the pairs $(k+1,m-k)$,
  $(k, m-k+1)$, $(k+2^q-1,m-k)$, and $(k,m-k+2^q-1)$
  have disjoint binary expansions.
  \end{itemize}

  \end{enumerate}

 \end{corol}

 \end{document}